\documentclass[aap,citesort,MSNbibl,dvips]{arximspdf}
\usepackage{mathbh}
\usepackage{graphicx}

\doi{10.1214/09-AAP630}
\volume{20}
\issue{3}
\pubyear{2010}
\firstpage{806}
\lastpage{840}

\makeatletter

\DeclareMathAlphabet\mathcaligr{OMS}{cmsy}{m}{n}

\newtheorem{proper}{Property}[section]
\newtheorem{lemm}[proper]{Lemma}
\newtheorem{theo}[proper]{Theorem}
\newtheorem{prop}[proper]{Proposition}
\newproclaim{defi}[proper]{Definition}

\makeatother

\begin{document}
\begin{frontmatter}

\title{Geography of local configurations}
\runtitle{Geography of local configurations}

\begin{aug}
\author[A]{\fnms{David} \snm{Coupier}\corref{}\ead[label=e1]{david.coupier@math.univ-lille1.fr}}
\runauthor{D. Coupier}
\affiliation{Universit\'{e} Lille 1}
\address[A]{Laboratoire Paul Painlev\'{e}\\
CNRS UMR 8524\\
Universit\'{e} Lille 1\\
59655 Villeneuve d'Ascq Cedex\\
France\\
\printead{e1}} 
\end{aug}

\received{\smonth{7} \syear{2007}}
\revised{\smonth{7} \syear{2009}}

%
\begin{abstract}
A $d$-dimensional binary Markov random field on a lattice torus is
considered. As the size $n$ of the lattice tends to infinity,
potentials $a=a(n)$ and $b=b(n)$ depend on $n$. Precise bounds for the
probability for local configurations to occur in a large ball are
given. Under some conditions bearing on $a(n)$ and $b(n)$, the distance
between copies of different local configurations is estimated according
to their weights. Finally, a sufficient condition ensuring that a given
local configuration occurs everywhere in the lattice is suggested.
\end{abstract}

%

\end{frontmatter}

\section{Introduction}
\label{section:intro}

\begin{figure}

\includegraphics{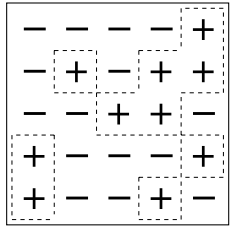}

\caption{A local configuration $\eta$ with
$k(\eta)=|V_{+}(\eta)|=10$ positive vertices and a perimeter
$\gamma(\eta)$ equals to $58$, in dimension $d=2$ and on a ball of
radius $r=2$ (with $\rho=1$ and relative to the $L_{\infty}$ norm).}
\label{fig:localconfig}
\end{figure}

In the theory of random graphs, inaugurated by Erd{\H{o}}s and
R\'{e}nyi \cite{ErdosRenyi},
the appearance of a given subgraph has been widely studied (see
Bollob\'{a}s \cite{Bollobas} and Spencer \cite{Spencer} for general references).
In the random graph formed by $n$ vertices, in which the edges are
chosen independently with pro\-bability $0<p<1$, a subgraph may occur
or not
according to the value of $p=p(n)$. In addition, under a certain
condition on the pro\-bability $p(n)$, its number of occurrences in the graph
is asymptotically (i.e., as $n\to+\infty$) Poissonian. Replacing the
edges with the states of a binary Markov random field, the notion of
subgraph corresponds to the notion of what we will call \textit{local
configuration}. Figure \ref{fig:localconfig} shows an example. Many
si\-tuations can be modeled by binary Markov random fields; a~vertex
and its state correspond to a pixel and its color (black or white) in
image analysis, to an individual and its opinion (yes or no) in
sociology or to an atom and its spin (positive or negative) in
statistical physics. This last interpretation leads to the well-known
Ising model. See \cite{KindermannSnell} for details.

%

Following the theory of random graphs, the appearance of a given local
configuration has been investigated in previous works (see \cite
{Coupier06-2} and \cite{CoupierDoukhanYcart}). In this article, this
study is extended into three directions. First, the speed at which
local configurations occur is specified. Moreover, when the number of
copies in the graph of a given local configuration is finite, the
states of vertices surrounding one of these copies are described.
Finally, a sufficient condition ensuring that a given local
configuration is present everywhere in the graph is stated. The results
obtained in these three directions are based on the same tools; the
Markovian character of the measure, the control of the conditional
probability for a local configuration to occur in the graph and the FKG
inequality \cite{FKG}.

Let us consider a lattice graph in dimension $d\geq1$, with periodic
boun\-dary conditions (lattice torus). The vertex set is $V_{n}=\{
0,\ldots,n-1\}^{d}$. The integer $n$ will be called the \textit{size}
of the lattice. The edge set, denoted by $E_{n}$, will be specified by
defining the set of neighbors $\mathcaligr{V}(x)$ of a given vertex $x$:
%
%
\begin{equation}
\label{def:voisinage}
\mathcaligr{V}(x) = \{ y \neq x \in V_{n} , \|y-x\|_{q} \leq\rho\} ,
\end{equation}
where the substraction is taken componentwise modulo $n$, $\|\cdot\|
_{q}$ stands for the $L_{q}$ norm in $\mathbb{R}^{d}$ ($1\leq q\leq
\infty$), and $\rho$ is a fixed integer. For instance, the square
lattice is obtained for $q=\rho=1$. Replacing the $L_{1}$ norm with
the $L_{\infty}$ norm adds the diagonals. From now on, all operations
on vertices will be understood modulo~$n$. In particular, each vertex
of the lattice has the same number of neighbors; we denote by
$\mathcaligr
{V}$ this number.

A \textit{configuration} is a mapping from the vertex set $V_{n}$ to
the state space $\{-1,+1\}$. Their set is denoted by $\mathcaligr
{X}_{n}=\{-1,+1\}^{V_{n}}$ and called the \textit{configuration set}.
In the following, we shall merely denote by $+$ and $-$ the states $+1$
and $-1$. Let $a$ and $b$ be two reals. The Gibbs measure associated to
potentials $a$ and $b$ is the probability measure $\mu_{a,b}$ on
$\mathcaligr{X}_{n}=\{-,+\}^{V_{n}}$ defined by: for all $\sigma\in
\mathcaligr{X}_{n}$,
%
%
\begin{equation}
\label{def:mesgibbs}
\mu_{a,b} (\sigma) = \frac{1}{Z_{a,b}} \exp\biggl( a \sum_{x \in V_{n}}
\sigma(x) + b \sum_{ \{x,y\} \in E_{n}} \sigma(x) \sigma(y)
\biggr) ,
\end{equation}
where the normalizing constant $Z_{a,b}$ is such that $\sum_{\sigma
\in\mathcaligr{X}_{n}} \mu_{a,b}(\sigma)=1$. Expectations relative to
$\mu_{a,b}$ will be denoted by $\mathbb{E}_{a,b}$. Georgii \cite{Georgii}
and Malyshev and Minlos \cite{Malyshev} constitute classical
references on Gibbs measures.

Throughout this paper, some hypotheses on $a$ and $b$ are made. The
model remaining unchanged by swapping positive and negative vertices
and replacing $a$ by $-a$, we chose to study only negative values of
the potential $a$. Thus, in order to use the FKG inequality, the
potential $b$ is supposed nonnegative. Finally, as the size $n$ of the
lattice tends to infinity, $a=a(n)$ and $b=b(n)$ are allowed to depend
on $n$. The case where $a(n)$ tends to $-\infty$ corresponds to rare
positive vertices among a majority of negative ones. So as to simplify
formulas, the Gibbs measure $\mu_{a(n),b(n)}$ is still denoted by $\mu
_{a,b}$.

In statistical physics, which is the point of view of \cite
{Coupier06-2}, the probabilistic model previously defined corresponds
to the ferromagnetic Ising model. In this context, potentials $a$ and
$b$ are, respectively, called the \textit{magnetic field} and the
\textit{pair potential}.

We are interested in the appearence in the graph $G_{n}$ of families of
local configurations. See Section \ref{section:localconfig} for a
precise definition and Figure \ref{fig:localconfig} for an example.
Such configurations are called ``local'' in the sense that the vertex
set on which they are defined is fixed and does not depend on $n$. A
local configuration $\eta$ is determined by its set of positive
vertices $V_{+}(\eta)$ whose cardinality and perimeter are,
respectively, denoted by $k(\eta)$ and $\gamma(\eta)$. A natural idea
(coming from \cite{Coupier06-2}) consists in regarding both parameters
$k(\eta)$ and $\gamma(\eta)$ through the same quantity; the \textit
{weight} of the local configuration $\eta$
\[
W_{n}(\eta) = \exp\bigl( 2 a(n) k(\eta) - 2 b(n) \gamma(\eta)
\bigr) .
\]
This notion plays a central role in our study. Indeed, the weight
$W_{n}(\eta)$ represents the probabilistic cost associated to a given
occurrence of $\eta$.

Proving some sharp inequalities is generally more difficult than
stating only limits. In the case of random graphs, Janson, {\L}uczak
and Ruci{\'n}ski \cite
{JLR}, thus Janson \cite{Janson}, have obtained exponential bounds for
the probability of nonexistence of subgraphs. Some other useful
inequalities have been suggested by Boppona and Spencer \cite
{Boppona}. In bond percolation on $\mathbb{Z}^{d}$, it is believed
that, in the subcritical phase, the probability for the radius of an
open cluster of being larger than $n$ behaves as an exponential term
multiplied by a power of $n$; see Grimmett \cite{Grimmett}, page 85, for
precise bounds. But, when the variables of the system are dependent, as
the states of a Markov random field, such inequalities become harder to
obtain. The Stein--Chen method (see Barbour, Holst and Janson \cite
{Barbour} for a
very complete reference or \cite{Chen} for the original paper of Chen)
is a useful way to bound the error of a Poisson approximation and so,
in particular, to bound the absolute value of the difference between
the probability of a property of the model and its limit. An example of
such a property is the appearance of a negative vertex (see Ganesh et
al. \cite{GHOSU}) or more generally that of any given local
configuration \cite{Coupier06-1}. These two previous papers concern
the case of a divergent potential $|a|$ and a constant potential $b$.
Coupling with the loss-network space-time representation due to
Fern{\'a}ndez, Ferrari and Garcia \cite{FFG2}, Ferrari and Picco \cite
{FerrariPicco}
have proved an exponential bound for large contours at low temperature
(i.e., $b$ large enough) and zero magnetic field (i.e., $a=0$).
Finally, under mixing conditions, various exponential approximations
with error bounds have been proved; see Abadi and Galves \cite
{AbadiGalves} for an overview and Abadi et al. \cite{ACRV} for the
high temperature case (i.e., $b$ small enough).

Our first goal is to establish precise lower and upper bounds for the
probability for certain families of local configurations to occur in
the graph for unbounded potentials $a(n)$ and $b(n)$.

Let $x\in V_{n}$ be a vertex, $W\in\,]0,1[$ be a real number and $R$ be
an integer. We denote by $\mathcaligr{A}(x,R,W)$ the (interpreted) event
``a local configuration whose weight is smaller than $W$ occurs
somewhere in the ball of center $x$ and radius $R$.''
The study of the probability of
this event becomes interesting when the radius $R=R(n)$ and the weight
$W=W_{n}$ depend on $n$ and tend, respectively, to infinity and zero.
Then, Theorem \ref{theo:expobounds} gives precise bounds for the
probability of the opposit event $^{c}\mathcaligr{A}(x,R(n),W_{n})$. There
exist two (explicit) constants $K>K'>0$ such that for $n$ large enough,
\[
\exp( -K R(n)^{d} W_{n} ) \leq\mu_{a,b} (
^{c}\mathcaligr{A}(x,R(n),W_{n}) ) \leq\exp( -K' R(n)^{d}
W_{n} ) .
\]
Another interesting problem consists in describing geographically what
happens in the studied model: the size of objects occurring in the
model and the distance between them. The size of components in random
graphs or the radius of open clusters in percolation are two classical
examples (see, respectively, \cite{Bollobas} and \cite{Grimmett}).
For the low-temperature plus-phase of the Ising model, Chazottes and
Redig \cite{ChazottesRedig} have studied the appearance of the first
two copies of a given pattern in terms of occurrence time and
repetition time. The occurrence time $T_{A}$ of a pattern $A$
represents the volume of the smallest set of vertices in which $A$ can
be found. As the size of the pattern increases, the distribution of
$T_{A}$ is approximated by an exponential law with error bounds. The
same is true for the repetition time $R_{A}$. Similar results exist for
sufficiently mixing Gibbs random fields (see \cite{ACRV}).

In our context, the studied objects are local. Hence, our second goal
is to estimate the distance between different copies of local
configurations occurring in the graph.

Let $\eta$ be a local configuration. It has been proved in \cite
{Coupier06-2} that the number of copies of $\eta$ occurring in $G_{n}$
is asymptotically Poissonian provided the pro\-duct $n^{d}W_{n}(\eta)$
is constant [the precise hypotheses are denoted by (\ref{hypH}) and recalled
in the beginning of Section \ref{section:distancebetween}]. In
particular, the number of copies of $\eta$ in the graph is finite with
probability tending to $1$. Let $\eta_{x}$ be one of them (that
occurring on the ball centered at $x$). We first observe that vertices
surrounding $\eta_{x}$ are all negative (Lemma \ref{lemm:autour}).
Hence, a natural question is the distance from $\eta_{x}$ to the
closest positive vertices. Theorem \ref{theo:dist1} answers to this
question; it states that the distance from $\eta_{x}$ to the closest
local configurations of weight $W_{n}$ is of order $W_{n}^{-1/d}$. Two
other results complete the study of the geography of local
configurations under the hypothesis (\ref{hypH}). Theorem \ref{theo:dist2}
claims the distance between the closest local configurations (to $\eta
_{x}$) of weight $W_{n}$ is also of order $W_{n}^{-1/d}$. The situation
described by Theorems \ref{theo:dist1} and \ref{theo:dist2} is
represented in Figure \ref{fig:dist}. Finally, Theorem \ref
{theo:dist1} implies the distance between any two two any copies of
$\eta$ should be of order $n$, which is the size of the graph.
Proposition \ref{prop:dist0(n)} precises this intuition.

From \cite{Coupier06-2}, a condition ensuring that a local
configuration $\eta$ occurs in $G_{n}$ is deduced; if the product
$n^{d}W_{n}(\eta)$ tends to infinity then, with probability tending to
$1$, at least one copy of $\eta$ can be found somewhere in the graph.
However, an uncertainty remains about the places in $G_{n}$ where $\eta
$ occurs. A richer information would be to know when the local
configuration $\eta$ occurs everywhere in $G_{n}$; we will talk about
\textit{ubiquity} of $\eta$. Inequalities stated in the proof of
Theorem \ref{theo:expobounds} allow us to obtain such an information.

For that purpose, the lattice $V_{n}$ is divided into blocks of
$(2R(n)+1)^{d}$ vertices. Thus, a supergraph $\tilde{G}_{n}$ whose set
of vertices $\tilde{V}_{n}$ is formed by the centers of these blocks
is defined. In order to study the appearance of $\eta$ in each of
these blocks, the set of configurations of $\tilde{V}_{n}$ is endowed
with an appropriate measure $\tilde{\mu}_{a,b}$ (depending on $\eta
$). Proposition \ref{prop:ubiquity} precises the asymptotic behavior
of $\tilde{\mu}_{a,b}$ according to the weight $W_{n}(\eta)$, the
radius of the blocks $R(n)$ and the size $n$. In particular, if
\[
\lim_{n \to+\infty} R(n)^{d} \ln\biggl( \frac{n}{R(n)} \biggr)
^{-1} W_{n}(\eta) = +\infty,
\]
then, with probability tending to $1$, all the blocks contain at least
one copy of the local configuration $\eta$.

The paper is organized as follows. The notion of local configuration
$\eta$ is
defined in Section \ref{section:localconfig}. Its number of positive vertices
$k(\eta)$, its perimeter $\gamma(\eta)$ and its weight $W_{n}(\eta
)$ are also
introduced. Section \ref{section:tools} is devoted to the three main
tools of our study. Property \ref{proper:markovian} underlines the
Markovian character of the Gibbs measure $\mu_{a,b}$. A control of the
conditional probability for a local configuration to occur on a ball
uniformly on the neighborhood of that ball is given in Lemma \ref
{lemm:controlproba}. Finally, the FKG inequality is discussed at the
end of Section \ref{section:tools}. Section \ref{section:expobounds}
gives the proof of Theorem \ref{theo:expobounds}. The geography of
local configurations occurring in the graph is described in Section
\ref{section:distancebetween}. Section \ref{subsection:motivation}
introduces the problem [the hypothesis (\ref{hypH}) and Lemma \ref
{lemm:autour}]. Theorems \ref{theo:dist1} and \ref{theo:dist2} and
Proposition \ref{prop:dist0(n)} are stated in Section~\ref
{subsection:main} and proved in Section \ref{subsection:proofs}. The
graph $\tilde{G}_{n}$ and the measure $\tilde{\mu}_{a,b}$ are
defined in Section \ref{section:ubiquity}. The latter inherits from
$\mu_{a,b}$ its Markovian character (Property \ref
{proper:markovian_tilde}). Finally, Proposition \ref{prop:ubiquity}
studies its asymptotic behavior.

\section{Local configurations}
\label{section:localconfig}


Let us start with some notation and definitions. Given
$\zeta\in\mathcaligr{X}_{n}=\{-,+\}^{V_{n}}$ and $V\subset V_{n}$, we
denote by
$\zeta_{V}$ the natural projection of $\zeta$ over $\{-,+\}^{V}$. If
$U$ and
$V$ are two disjoint subsets of $V_{n}$ then $\zeta_{U}\zeta'_{V}$ is the
configuration on $U\cup V$ which is equal to $\zeta$ on $U$ and $\zeta
'$ on
$V$. Let us denote by $\delta V$ the neighborhood of $V$ [corresponding to
(\ref{def:voisinage})]:
\[
\delta V = \bigl\{ y \in V_{n} \setminus V , \exists x \in V , \{ x,y
\} \in
E_{n} \bigr\}
\]
and by $\overline{V}$ the union of the two disjoint sets $V$ and
$\delta
V$. Moreover, $|V|$ denotes the cardinality of $V$ and $\mathcaligr
{F}(V)$ the
$\sigma$-algebra generated by the configurations of $\{-,+\}^{V}$.
Finally, if $A\in\mathcaligr{F}(V_{n})$, we denote by $^{c}A$ the
opposit event.

As usual, the graph distance $\mathrm{dist}$ is defined as the minimal
length of
a path
between two vertices. We shall denote by $B(x,r)$ the ball of center
$x$ and
radius $r$:
\[
B(x,r) = \{ y\in V_{n} ; \operatorname{dist}(x,y)\leq r \} .
\]
In the case of balls, $\overline{B(x,r)}=B(x,r+1)$. In order to avoid
unpleasant situations, like self-overlapping balls, we will always
assume that
$n>2\rho r$. If $n$ and $n'$ are both larger than $2\rho r$, the balls
$B(x,r)$ in $G_{n}$ and $G_{n'}$ are isomorphic. Two properties of the balls
$B(x,r)$ will be crucial in what follows. The first one is that two
balls with
the same radius are translates of each other:
\[
B(x+y,r)=y+B(x,r).
\]
The second one is that for $n>2\rho r$, the cardinality of $B(x,r)$ depends
only on $r$ and neither on $x$ nor on $n$: it will be denoted by
$\beta_{r}$. Observe that whatever the choices of $q$ and $\rho$ in
(\ref{def:voisinage}) the ball $B(x,r)$ is included in the sublattice
$[x-\rho r, x+\rho r]$, and thus
\[
\beta_{r} \leq( 2 \rho r + 1 )^{d} .
\]


Let $r$ be a positive integer, and consider a fixed ball with radius
$r$, say
$B(0,r)$. We denote by $\mathcaligr{C}_{r}=\{-,+\}^{B(0,r)}$ the set of
configurations on that ball. Elements of $\mathcaligr{C}_{r}$ will be called
\textit{local configurations with radius $r$}, or merely \textit{local
configurations} whenever the radius $r$ will be fixed. Of course, there
exists only a finite number of such configurations (precisely
$2^{\beta_{r}}$). See Figure \ref{fig:localconfig} for an example. Throughout
this paper, the radius $r$ will be constant, that is, it will not
depend on the
size $n$. Hence, defining local configurations on balls of radius $r$ serves
only to ensure that studied objects are ``local.'' In what follows,
$\eta$,
$\eta'$ will denote local configurations of radius~$r$.

A local configuration $\eta\in\mathcaligr{C}_{r}$ is determined by its subset
$V_{+}(\eta)\subset B(0,r)$ of positive vertices:
\[
V_{+}(\eta) = \{ x \in B(0,r), \eta(x) = + \} .
\]
The cardinality of this set will be denoted by $k(\eta)$ and its
complement in
$B(0,r)$, that is, the set of negative vertices of $\eta$, by
$V_{-}(\eta)$. Moreover, the geometry (in the sense of the graph
structure) of
the set $V_{+}(\eta)$ needs to be described. Let us define the
\textit{perimeter} $\gamma(\eta)$ of the local configuration $\eta$
by the
formula
\[
\gamma(\eta) = \mathcaligr{V} | V_{+}(\eta) | - 2 \bigl| \bigl\{ \{
x,y\} \in
V_{+}(\eta) \times V_{+}(\eta), \{x,y\} \in E_{n} \bigr\} \bigr| ,
\]
where $\mathcaligr{V}$ is the number of neighbors of a vertex. In other words,
$\gamma(\eta)$ counts the pairs of neighboring vertices $x$ and $y$ of
$B(0,r)$ having opposite spins (under $\eta$) and those such that
$x\in
B(0,r)$, $y\in\delta B(0,r)$ and $\eta(x)=+$.

We denote by $W_{n}(\eta)$ and call the \textit{weight} of the local
configuration $\eta$ the following quantity:
\[
W_{n}(\eta) = \exp\bigl( 2 a(n) k(\eta) - 2 b(n) \gamma(\eta)
\bigr) .
\]
Since $a(n)<0$ and $b(n)\geq0$, the weight $W_{n}(\eta)$ satisfies
$0<W_{n}(\eta)\leq1$. That of the local configuration having only negative
vertices, denoted by $\eta_{-}$ and called the \textit{negative local
configuration}, is equal to $1$. If $\eta\not=\eta_{-}$ then $k(\eta
)\geq1$
and $\gamma(\eta)\geq\mathcaligr{V}$. It follows that
\[
W_{n}(\eta) \leq\exp\bigl( 2 a(n) - 2 b(n) \mathcaligr{V} \bigr) .
\]
Actually, the weight $W_{n}(\eta)$ represents the probabilistic cost
associated to the presence of $\eta$ on a given ball. This idea will be
clarified in the next section (Lem\-ma~\ref{lemm:controlproba}).

Remark the notation $k(\cdot)$, $\gamma(\cdot)$ and $W_{n}(\cdot)$
can be
naturally extended to any configuration $\zeta\in\{-,+\}^{V}$,
$V\subset
V_{n}$.

Let $\eta\in\mathcaligr{C}_{r}$. For each vertex $x\in V_{n}$, denote by
$\eta_{x}$ the translation of $\eta$ onto the ball $B(x,r)$ (up to periodic
boundary conditions):
\[
\forall y \in V_{n}\qquad \operatorname{dist}(0,y)\leq r \quad
\Longrightarrow\quad\eta_{x}(x+y) =
\eta(y) .
\]
In particular, $V_{+}(\eta_{x})=x+V_{+}(\eta)$. So, $\eta$ and $\eta
_{x}$ have
the same number of positive vertices and the same perimeter. So do their
weights. Let us denote by $I_{x}^{\eta}$ the indicator function
defined on
$\mathcaligr{X}_{n}$ as follows: $I_{x}^{\eta}(\sigma)$ is $1$ if the
restriction
of the configuration $\sigma\in\mathcaligr{X}_{n}$ to the ball $B(x,r)$ is
$\eta_{x}$ and $0$ otherwise.


Let $V$ and $V'$ be two disjoint subsets of vertices. The following relations
\[
k(\zeta\zeta') = k(\zeta) + k(\zeta') \quad\mbox{and}\quad
\gamma(\zeta\zeta') \leq\gamma(\zeta) + \gamma(\zeta')
\]
are true whatever the configurations $\zeta\in\{-,+\}^{V}$ and
$\zeta'\in\{-,+\}^{V'}$. As an immediate consequence, the weight
$W_{n}(\zeta
\zeta')$ is larger than the product $W_{n}(\zeta)W_{n}(\zeta')$. The
\textit{connection} between $\zeta\in\{-,+\}^{V}$ and $\zeta'\in\{
-,+\}^{V'}$,
denoted by $\operatorname{conn}(\zeta,\zeta')$, is defined by
\[
\operatorname{conn}(\zeta,\zeta') = \bigl| \bigl\{ \{y,z\} \in E_{n} ,
y \in V , z \in V'
\mbox{ and } \zeta(y) = \zeta'(z) = + \bigr\} \bigr| .
\]
This quantity allows us to link the perimeters of the configurations
$\zeta\zeta'$, $\zeta$ and $\zeta'$ together:
\[
\gamma(\zeta\zeta') + 2 \operatorname{conn}(\zeta,\zeta') = \gamma(\zeta
) +
\gamma(\zeta')
\]
and therefore their weights:
%
%
\begin{equation}
\label{lienweights}
W_{n}(\zeta\zeta') \exp(- 4 b(n) \operatorname{conn}(\zeta,\zeta')) =
W_{n}(\zeta)
W_{n}(\zeta') .
\end{equation}
In particular, if the connection $\operatorname{conn}(\zeta,\zeta')$ is
null then the
weight $W_{n}(\zeta\zeta')$ is equal to the product
$W_{n}(\zeta)W_{n}(\zeta')$. This is the case when $\overline{V}\cap V'
=\varnothing$.

\section{The three main tools}
\label{section:tools}

This section is devoted to the main tools on which are based all the results
of this paper: the Markovian character of the Gibbs measure $\mu
_{a,b}$, a
control of the probability for $\eta\in\mathcaligr{C}_{r}$ to occur on
a given
ball and the FKG inequality. Except the first part of Lemma
\ref{lemm:controlproba}, the results of this section are already
known.


Two subsets of vertices $U$ and $V$ of $V_{n}$ are said
$\mathcaligr{V}$-\textit{disjoint} if none of the vertices of $U$ belong
to the
neighborhood of one of the vertices of $V$. In other words, $U$ and $V$ are
$\mathcaligr{V}$-disjoint if and only if $U\cap\overline{V}=\varnothing
$ (or
equivalently $\overline{U}\cap V=\varnothing$). For example, two balls $B(x,r)$
and $B(x',r)$ are $\mathcaligr{V}$-disjoint if and only if the distance between
their centers $x$ and $x'$ is larger than $2r+1$, that is,
$\operatorname{dist}(x,x')>2r+1$.

The following result is a classical property of Gibbs measures (see
\cite{Georgii}, page~157); it describes the Markovian character of
$\mu_{a,b}$. The second part of Property \ref{proper:markovian} means that,
given two $\mathcaligr{V}$-disjoint sets $U$ and $V$, the $\sigma$-algebras
$\mathcaligr{F}(U)$ and $\mathcaligr{F}(V)$ are conditionally
independent knowing
the configuration on $\delta U\cup\delta V$.

For any sets of vertices $V,V'$ and for any event $A\in\mathcaligr
{F}(V)$, the
function $\mu_{a,b}(A|\mathcaligr{F}(V'))$ denotes the
$\mathcaligr{F}(V')$-measurable random variable defined as follows; for
$\sigma\in\{-,+\}^{V'}$, $\mu_{a,b}(A|\mathcaligr{F}(V'))(\sigma)$ is the
conditional probability $\mu_{a,b}(A|\sigma)$.
\begin{proper}
\label{proper:markovian}
Let $V,V'\subset V_{n}$ be two sets of vertices such that $V\cap
V'=\varnothing$ and $\delta V\subset V'$. Then, for all $A\in\mathcaligr{F}(V)$,
%
%
\begin{equation}
\label{markovian1}
\mu_{a,b}(A | \mathcaligr{F}(V')) = \mu_{a,b}(A | \mathcaligr{F}(\delta
V)) .
\end{equation}
Let $U,V,V'\subset V_{n}$ be three sets of vertices such that $U$ and
$V$ are $\mathcaligr{V}$-disjoint, $(U\cup V)\cap V'=\varnothing$ and
$\delta U\cup\delta V\subset V'$. Then, for all $A\in\mathcaligr{F}(U)$
and $B\in\mathcaligr{F}(V)$,
%
%
\begin{equation}
\label{markovian2}
\mu_{a,b}\bigl(A \cap B | \mathcaligr{F}(V')\bigr) = \mu_{a,b}(A |
\mathcaligr
{F}(V'))
\mu_{a,b}(B | \mathcaligr{F}(V')) .
\end{equation}
\end{proper}

Let us note that (\ref{markovian1}) is a consequence of the identity
\[
\forall\sigma\in\{-,+\}^{\delta V} , \forall\sigma' \in\{-,+\}^{V'
\setminus\delta V} \qquad \mu_{a,b}(A | \sigma\sigma') = \mu
_{a,b}(A |
\sigma) ,
\]
which itself relies on the exponential form of the Gibbs measure $\mu_{a,b}$
[see (\ref{def:mesgibbs})]. The proof of a similar identity is available
in \cite{Malyshev}, page 7.

Besides, the second part of Property \ref{proper:markovian} can be
immediatly extended to any finite family of sets of vertices which are
two by two
$\mathcaligr{V}$-disjoint. This remark will be used often in the
following.


Let $\eta$ be a local configuration. Thanks to the translation
invariance of the graph~$G_{n}$, the indicator functions $I_{x}^{\eta
}$, $x\in V_{n}$, have the same distribution. So, let us pick a vertex
$x$. For any configuration $\sigma\in\{-,+\}^{\delta B(x,r)}$, the
quantity $\mu_{a,b}(I_{x}^{\eta}=1|\sigma)$ represents the
probability for $\eta$ (or $\eta_{x}$) to occur on the ball $B(x,r)$
knowing what happens on its neighborhood. A precise study of this
conditional probability has been done in \cite{Coupier06-2}, Section
2. From there, a control is deduced:
\begin{lemm}
\label{lemm:controlproba}
Let $\eta$ be a local configuration with radius $r$, and $x$ be a vertex.

On the one hand, if $a(n)+\mathcaligr{V}b(n)\leq0$ then there exists a constant
$c_{r}>0$ such that, for all configuration $\sigma\in\{-,+\}^{\delta B(x,r)}$,
%
%
\begin{equation}
\label{controlinfproba}
\mu_{a,b}(I_{x}^{\eta} = 1 | \sigma) \geq c_{r} W_{n}(\eta) .
\end{equation}
In particular, $\mu_{a,b}(I_{x}^{\eta}=1)$ satisfies the same lower
bound. On
the other hand, if $a(n)+2\mathcaligr{V}b(n)\leq0$ then there exists a constant
$C_{r}>0$ such that
%
%
\begin{equation}
\label{controlsupproba}
\mu_{a,b}(I_{x}^{\eta} = 1) \leq C_{r} W_{n}(\eta) .
\end{equation}
\end{lemm}

The constants $C_{r}$ and $c_{r}$ depend on the radius $r$ and on
parameters $q$, $\rho$ and $d$ but not on the local configuration
$\eta$ nor on $n$.
\begin{pf*}{Proof of Lemma \protect\ref{lemm:controlproba}}
Let $\eta$ be a local configuration with radius $r$, $x$ be a vertex
and $\sigma$ be an element of $\{-,+\}^{\delta B(x,r)}$. Lemma 2.2 of
\cite{Coupier06-2} allows us to write
the conditional probability $\mu_{a,b}(I_{x}^{\eta}=1|\sigma)$ as a
function of weights of some configurations:
%
%
\begin{eqnarray}
\label{withweights}
\mu_{a,b}(I_{x}^{\eta} = 1 | \sigma) & = & \frac{W_{n}(\eta_{x}
\sigma)}{\sum_{\eta'\in\mathcaligr{C}_{r}} W_{n}(\eta'_{x} \sigma)}
\nonumber\\[-8pt]\\[-8pt]
& \geq& W_{n}(\eta) \biggl( \sum_{\eta'\in\mathcaligr{C}_{r}}
\frac{W_{n}(\eta'_{x} \sigma)}{W_{n}(\sigma)} \biggr)^{-1} ,\nonumber
\end{eqnarray}
using the inequality $W_{n}(\eta_{x}\sigma)\geq W_{n}(\eta
)W_{n}(\sigma)$. Let $\eta'$ be a local configuration. Its perimeter
$\gamma(\eta')$ can be bounded as follows:
\[
\mathcaligr{V} k(\eta') \geq\gamma(\eta') \geq\operatorname{conn}(\eta
'_{x},\sigma
) .
\]
Hence, $\gamma(\eta')+\mathcaligr{V}k(\eta')$ is larger than
$2\operatorname{conn}(\eta'_{x},\sigma)$ and the difference $\gamma(\eta
'_{x}\sigma)-\gamma(\sigma)$ satisfies
\begin{eqnarray*}
\gamma(\eta'_{x}\sigma)-\gamma(\sigma) & = & \gamma(\eta') - 2
\operatorname{conn}(\eta'_{x},\sigma) \\
& \geq& - \mathcaligr{V} k(\eta') .
\end{eqnarray*}
We deduce from this last inequality that
\begin{eqnarray*}
\frac{W_{n}(\eta'_{x} \sigma)}{W_{n}(\sigma)} & = & \exp\bigl(
2a(n)k(\eta')
- 2b(n)\bigl(\gamma(\eta'_{x}\sigma) - \gamma(\sigma)\bigr) \bigr) \\
& \leq& \exp\bigl( \bigl(2a(n) + 2\mathcaligr{V}b(n)\bigr) k(\eta')
\bigr) .
\end{eqnarray*}
Now, the hypothesis $a(n)+\mathcaligr{V}b(n)\leq0$ implies the ratio
$W_{n}(\eta'_{x}\sigma)$ divided by $W_{n}(\sigma)$ is smaller than
$1$. Thanks to (\ref{withweights}), the conditional probability $\mu
_{a,b}(I_{x}^{\eta}=1|\sigma)$ is larger than $|\mathcaligr
{C}_{r}|^{-1}W_{n}(\eta)$; $c_{r}=|\mathcaligr{C}_{r}|^{-1}=2^{-\beta
_{r}}$ is suitable.

Since the lower bound of (\ref{controlinfproba}) is uniform on the
configuration $\sigma\in\{-,+\}^{\delta B(x,r)}$, the same inequality
holds for
\[
\mu_{a,b}(I_{x}^{\eta} = 1) = \mathbb{E}_{a,b} \bigl[\mu_{a,b}
\bigl(
I_{x}^{\eta}
= 1 | \mathcaligr{F}(\delta B(x,r)) \bigr) \bigr] .
\]
Finally, (\ref{controlsupproba}) has been proved in Proposition 3.2 of
\cite{Coupier06-2} with $C_{r}=2^{\mathcaligr{V}\beta_{r}}$.
\end{pf*}

The lower bound given by (\ref{controlinfproba}) has the advantage of being
uniform on the configuration $\sigma\in\{-,+\}^{\delta B(x,r)}$. In addition,
there is no uniform upper bound for the conditional probability
$\mu_{a,b}(I_{x}^{\eta}=1|\sigma)$. In order to make up for this
gap, we will
have recourse to the FKG inequality.


There is a natural partial ordering on the configuration set
$\mathcaligr{X}_{n}=\{-,+\}^{V_{n}}$ defined by $\sigma\leq\sigma'$ if
$\sigma(x)\leq\sigma'(x)$ for all vertices $x\in V_{n}$. A real function
defined on $\mathcaligr{X}_{n}$ is \textit{increasing} if $f(\sigma
)\leq
f(\sigma')$ whenever $\sigma\leq\sigma'$. An event
$A\in\mathcaligr{F}(V_{n})$ is also said \textit{increasing} whenever
its indicator function $f=\mathbh{1}_{A}$ is increasing. Conversely, a
\textit{decreasing} event is an event whose complementary set [in
$\mathcaligr{F}(V_{n})$] is increasing.

Let us focus on an example of increasing event which will be central in our
study. For $0\leq W\leq1$, we denote by $\mathcaligr{C}_{r}(W)$ the set
of local
configurations with radius $r$ whose weight is smaller than $W$:
\[
\mathcaligr{C}_{r}(W) = \{ \eta\in\mathcaligr{C}_{r} , W_{n}(\eta)
\leq W \} .
\]
Thus, let us consider two local configurations $\eta$, $\eta'$ such
that the
set of positive vertices of $\eta'$ contains that of $\eta$. The perimeter
$\gamma(\eta')$ is not necessary larger than $\gamma(\eta)$: roughly,
$V_{+}(\eta)$ may have holes. However, the inequality
\[
\gamma(\eta') \geq\gamma(\eta) - \bigl(k(\eta') - k(\eta)\bigr)
\mathcaligr{V}
\]
holds. Hence, the ratio
\begin{eqnarray*}
\frac{W_{n}(\eta')}{W_{n}(\eta)} & = & \exp\bigl( 2 a(n) \bigl(k(\eta
') - k(\eta)\bigr)
- 2 b(n) \bigl(\gamma(\eta') - \gamma(\eta)\bigr) \bigr) \\
& \leq& \exp\bigl( \bigl( 2 a(n) + 2 \mathcaligr{V} b(n) \bigr)
\bigl(k(\eta') -
k(\eta)\bigr) \bigr)
\end{eqnarray*}
is smaller than $1$ whenever $a(n)+\mathcaligr{V}b(n)$ is negative. As a
consequence, under this hypothesis, the set $\mathcaligr{C}_{r}(W)$
allows us to
build some increasing events; for instance,
\[
\bigcup_{\eta\in\mathcaligr{C}_{r}(W)} \{ I_{x}^{\eta} = 1
\} .
\]
For a positive value of the pair potential $b$, the Gibbs measure $\mu_{a,b}$
defined by (\ref{def:mesgibbs}) satisfies the FKG inequality, that is,
%
%
\begin{equation}
\label{inegFKG}
\mu_{a,b} (A\cap A') \geq\mu_{a,b} (A) \mu_{a,b} (A')
\end{equation}
for all increasing events $A$ and $A'$. See, for instance, Section 3 of
\cite{FKG}. In statistical physics, the hypothesis $b\geq0$
corresponds to \textit{ferromagnetic interaction}. Note
that inequality (\ref{inegFKG}) can be easily extended to decreasing
events (see \cite{Grimmett}). Moreover, the union and the intersection
of increasing events are still increasing. The same is true for
decreasing events. Hence, the FKG inequality applies to any family made
up of a finite number of decreasing events $A_{1},\ldots,A_{m}$:
\[
\mu_{a,b} (A_{1} \cap\cdots\cap A_{m}) \geq\mu_{a,b} (A_{1})
\cdots\mu_{a,b} (A_{m}) .
\]

\section{Exponential bounds for the probability of nonexistence}
\label{section:expobounds}

Let $x\in V_{n}$ be a vertex, $R\geq r\geq1$ be two integers and $W$
be a
positive real. Let us denote by $\mathcaligr{A}(x,R,W)$ the following event:
\[
\exists y \in B(x,R-r) , \exists\eta\in\mathcaligr{C}_{r}(W) \qquad
I_{y}^{\eta} = 1 .
\]
The event $\mathcaligr{A}(x,R,W)$ means at least one copy of a local
configuration (with radius $r$) whose weight is smaller than $W$ can be found
somewhere in the large ball $B(x,R)$. Theorem \ref{theo:expobounds} gives
exponential bounds for the probability of the opposite event
$^{c}\mathcaligr{A}(x,R(n),W_{n})$.
\begin{theo}
\label{theo:expobounds}
Assume that the magnetic field $a(n)$ is negative, the pair potential $b(n)$
is nonnegative and they satisfy $a(n)+2\mathcaligr{V}b(n)\leq0$. Let
$(W_{n})_{n\in\mathbb{N}}$ be a sequence of positive reals satisfying the
following property: there exist an integer $N$ and $0<\epsilon<1$ such that
\[
\forall n \geq N \qquad \min_{\eta\in\mathcaligr{C}_{r}} W_{n}(\eta)
\leq W_{n}
\leq\epsilon c_{r} C_{r}^{-1} .
\]
Let $K=(2\rho)^{d}(1-\epsilon)^{-1}c_{r}^{-1}C_{r}$ and $K'=\tau
c_{r}$. Then, for all
$n\geq N$ and for all vertex $x$,
%
%
\begin{equation}\quad
\label{expobounds}
\exp( -K R(n)^{d} W_{n} ) \leq\mu_{a,b} (
^{c}\mathcaligr{A}(x,R(n),W_{n}) ) \leq\exp( -K' R(n)^{d}
W_{n} ) .
\end{equation}
\end{theo}

Remark the constants $K$ and $K'$ only depend on $r$, $\epsilon$ and
parameters $q$, $\rho$ and $d$ (the constant $\tau$ will be
introduced in the
proof).

The fact that, for $n$ large enough, $W_{n}$ is assumed larger than the
smallest weight $W_{n}(\eta)$, $\eta\in\mathcaligr{C}_{r}$, only
serves to ensure that the set $\mathcaligr{C}_{r}(W_{n})$
[and the event $\mathcaligr{A}(x,R(n),W_{n})$ too] is nonempty.
Moreover, let us note hypothesis $W_{n}\leq\epsilon c_{r}C_{r}^{-1}$
is only used in the proof of the lower bound of (\ref{expobounds}).

The inequalities of (\ref{expobounds}) give a limit for the
probability of
$\mathcaligr{A}(x,R(n),W_{n})$:
\begin{eqnarray*}
&&\mbox{if } \lim_{n \to+\infty} R(n)^{d} W_{n} = 0 \qquad
\mbox{then } \lim_{n \to+\infty} \mu_{a,b} (\mathcaligr
{A}(x,R(n),W_{n})) = 0 ;
\\
&&\mbox{if } \lim_{n \to+\infty} R(n)^{d} W_{n} = +\infty
\qquad\mbox{then }
\lim_{n \to+\infty} \mu_{a,b} (\mathcaligr{A}(x,R(n),W_{n})) = 1 .
\end{eqnarray*}
Roughly speaking, if the radius $R(n)$ is small compared to
$W_{n}^{-1/d}$,
then asymptotically, with probability tending to $1$, there is no local
configuration whose weight is smaller than $W_{n}$. Conversely, if
$R(n)$ is
large compared to $W_{n}^{-1/d}$ then at least one copy of an element of
$\mathcaligr{C}_{r}(W_{n})$ can be found somewhere in $B(x,R(n))$, with
probability tending to $1$.

Finally, Theorem \ref{theo:expobounds} implies the quantity
%
%
\begin{equation}
\label{logproba}
\frac{1}{R(n)^{d} W_{n}} \ln( \mu_{a,b} (^{c}\mathcaligr
{A}(x,R(n),W_{n})) )
\end{equation}
belongs to the interval $[-K,-K']$. A natural question is to wonder whether
(\ref{logproba}) admits a limit as $n$ goes to infinity. It seems to
be difficult to answer to this question in general. However, in the
following particular case,
the answer is positive and the corresponding limit is known; see
\cite{Coupierthese} or \cite{Coupier06-1} for more details. Assume
the weight
$W_{n}=W_{n}(\eta)$ is the one of a local configuration $\eta$, the pair
potential $b(n)=b$ is a positive real number and the magnetic field $a(n)$
tends to $-\infty$. Under these conditions, the limit of the
probability for a
given local configuration to occur on the ball $B(x,R(n))$ essentially depends
on its number of positive vertices. Moreover, assume the radius $R(n)$
is such
that the product $R(n)^{d}W_{n}=R(n)^{d}W_{n}(\eta)$ is constant.
Then, the
quantity (\ref{logproba}) converges, as $n\to+\infty$, to
\[
- \mathop{\mathop{\sum}_{\eta' \in\mathcaligr{C}_{r}}}_
{k(\eta') =
k(\eta)} \exp\bigl( -2 b \bigl(\gamma(\eta') - \gamma(\eta)\bigr) \bigr).
\]
This section ends with the proof of Theorem \ref{theo:expobounds}.
\begin{pf*}{Proof of Theorem \protect\ref{theo:expobounds}}
Throughout this proof, potentials $a(n)$ and $b(n)$, and the sequence
$(W_{n})_{n\in\mathbb{N}}$ satisfy the hypotheses of Theorem
\ref{theo:expobounds}.

The two following remarks will lighten notation and formulas of the
proof. Thanks to the invariance translation of the graph $G_{n}$ it suffices
to prove Theorem \ref{theo:expobounds} for $x=0$. Hence, we merely
denote by
$\mathcaligr{A}(R(n),W_{n})$ the event $\mathcaligr{A}(0,R(n),W_{n})$.
Furthermore,
$\mathcaligr{A}(R(n),W_{n})$ and $\mathcaligr{A}(R(n),W_{n}')$ are
equal for
\[
W_{n}' = \max_{\eta\in\mathcaligr{C}_{r}(W_{n})} W_{n}(\eta) .
\]
As a consequence and without loss of generality, we can assume that,
for each
$n$, the weight $W_{n}$ belongs to $\{ W_{n}(\eta),
\eta\in\mathcaligr{C}_{r}\}$.


The proof of the lower bound of (\ref{expobounds}) requires the ferromagnetic
character of the probabilistic model, that is, the positivity of the
pair potential
$b(n)$. Since $a(n)+\mathcaligr{V}b(n)$ is negative, the event
\[
\bigcap_{\eta\in\mathcaligr{C}_{r}(W_{n})} \{ I_{y}^{\eta} = 0
\}
\]
is decreasing, whatever the vertex $y$. So the FKG inequality implies
%
%
\begin{equation}
\label{minorFKG}
\mu_{a,b} ( ^{c}\mathcaligr{A}(R(n),W_{n}) ) \geq\prod
_{y\in
B(0,R(n)-r)} \mu_{a,b} \biggl( \bigcap_{\eta\in\mathcaligr{C}_{r}(W_{n})}
\{ I_{y}^{\eta} = 0 \} \biggr) .
\end{equation}
Now, it suffices to control each term of the above product. This is the role
of Lem\-ma \ref{lemm:controlproba}. Let us pick $y\in B(0,R(n)-r)$. We get
\begin{eqnarray*}
\mu_{a,b} \biggl( \bigcap_{\eta\in\mathcaligr{C}_{r}(W_{n})}
\{ I_{y}^{\eta} = 0 \} \biggr) & = & 1 -
\sum_{\eta\in\mathcaligr{C}_{r}(W_{n})} \mu_{a,b} ( I_{y}^{\eta
} = 1 )
\\
& \geq& 1 - \sum_{\eta\in\mathcaligr{C}_{r}(W_{n})} C_{r} W_{n}(\eta
) \\
& \geq& 1 - c_{r}^{-1} C_{r} W_{n} ,
\end{eqnarray*}
since $|\mathcaligr{C}_{r}(W_{n})|\leq|\mathcaligr{C}_{r}|=c_{r}^{-1}$.
Let $N$ be
the integer introduced in the statement of Theorem \ref
{theo:expobounds} and
let $n\geq N$. Therefore, from (\ref{minorFKG}) and the inequality
$\ln(1-X)\geq-(1-\epsilon)^{-1}X$ valid for $0\leq X\leq\epsilon$,
it follows
that
\begin{eqnarray*}
\mu_{a,b} ( ^{c}\mathcaligr{A}(R(n),W_{n}) ) & \geq&
( 1 -
c_{r}^{-1} C_{r} W_{n} ) ^{(2\rho R(n))^{d}} \\
& \geq& \exp\bigl( -(2\rho)^{d}(1-\epsilon)^{-1} c_{r}^{-1} C_{r}
R(n)^{d} W_{n} \bigr).
\end{eqnarray*}
So, the lower bound of (\ref{expobounds}) is proved with
$K=(2\rho)^{d}(1-\epsilon)^{-1}c_{r}^{-1}C_{r}$.


In order to prove the upper bound of (\ref{expobounds}), let us denote by
$T_{n}$ the subset of $V_{n}$ defined by
\[
T_{n} = \biggl\{ i \bigl(\rho(2r+1)+ 1\bigr) , |i| = 0,\ldots,
\biggl\lfloor
\frac{\rho R(n)}{\rho(2r+1)+ 1} \biggr\rfloor- 1 \biggr\}^{d}
\cap B\bigl(0,R(n)-r\bigr) .
\]
Its cardinality, denoted by $\tau_{n}$, satisfies $\tau_{n}\geq\tau
R(n)^{d}$, for a positive constant $\tau>0$ not depending on the size
$n$. Let
$\mathcaligr{T}_{n}$ be the union of the balls with radius $r$ centered
at the
elements of $T_{n}$ (see Figure \ref{fig:reseau})
\[
\mathcaligr{T}_{n} = \bigcup_{y \in T_{n}} B(y,r)
\]
and denote by $\Gamma_{n}$ the event which means an element of
$\mathcaligr{C}_{r}(W_{n})$ occurs on a ball $B(y,r)$, $y\in T_{n}$:
\[
\Gamma_{n} = \bigcup_{y \in T_{n}} \bigcup_{\eta\in\mathcaligr
{C}_{r}(W_{n})} \{ I_{y}^{\eta} = 1 \} .
\]

%
\begin{figure}

\includegraphics{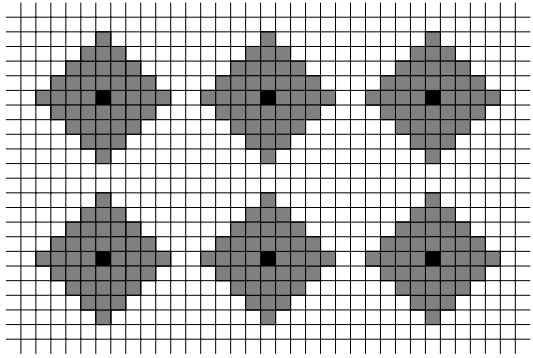}

\caption{The colored vertices represent the set
$\mathcaligr{T}_{n}$ in dimension $d=2$, with $r=2$, $\rho=2$ and
relative to
the $L_{1}$ norm. The black vertices are elements of $T_{n}$.}
\label{fig:reseau}
\end{figure}

Since the balls $B(y,r)$, $y\in T_{n}$, are included in $B(0,R(n))$,
the event
$\Gamma_{n}$ implies $\mathcaligr{A}(R(n),W_{n})$. So, the upper bound of
(\ref{expobounds}) follows from the next statement. Let $K'=\tau
c_{r}$. Then, for all $n$,
%
%
\begin{equation}
\label{stepupperbound}
\mu_{a,b} (^{c}\Gamma_{n}) \leq\exp( -K' R(n)^{d} W_{n}
) .
\end{equation}
Let us prove this inequality. Remark that any two balls $B(y,r)$ and
$B(y',r)$ where $y$ and $y'$ are distinct vertices of $T_{n}$ are
$\mathcaligr{V}$-disjoint. Indeed, by construction of the set $T_{n}$,
$\operatorname{dist}(y,y')>2r+1$. Hence, Property \ref
{proper:markovian} produces the
following identities:
%
%
\begin{eqnarray}
\label{utilisationmarkovian}
\mu_{a,b} ( ^{c}\Gamma_{n} ) & = & \mathbb{E}_{a,b}
[\mu_{a,b}
(^{c}\Gamma_{n} | \mathcaligr{F} ( \delta
\mathcaligr{T}_{n}
) ) ]\nonumber\\
& = & \mathbb{E}_{a,b} \biggl[\mu_{a,b} \biggl(
\bigcap_{y
\in T_{n}}
\bigcap_{\eta\in\mathcaligr{C}_{r}(W_{n})} I_{y}^{\eta} = 0 \Big|
\mathcaligr{F} ( \delta\mathcaligr{T}_{n} ) \biggr)
\biggr]
\nonumber\\[-8pt]\\[-8pt]
& = & \mathbb{E}_{a,b} \biggl[\prod_{y \in T_{n}} \mu_{a,b}
\biggl(
\bigcap_{\eta\in\mathcaligr{C}_{r}(W_{n})} I_{y}^{\eta} = 0
\Big| \mathcaligr{F} ( \delta\mathcaligr{T}_{n} ) \biggr)
\biggr]\nonumber\\
& = & \mathbb{E}_{a,b} \biggl[\prod_{y \in T_{n}} \mu_{a,b}
\biggl(
\bigcap_{\eta\in\mathcaligr{C}_{r}(W_{n})} I_{y}^{\eta} = 0
\Big| \mathcaligr{F} ( \delta B(y,r) ) \biggr)
\biggr] .\nonumber
\end{eqnarray}
Now, let $y\in T_{n}$ be a vertex and $\sigma\in\{-,+\}^{\delta
B(y,r)}$ be a
configuration. We can write
%
%
\begin{eqnarray}
\label{unmoinsprobacond}
\mu_{a,b} \biggl( \bigcap_{\eta\in\mathcaligr{C}_{r}(W_{n})}
I_{y}^{\eta}
= 0 \Big| \sigma\biggr) & = & 1 - \mu_{a,b} \biggl(
\bigcup_{\eta\in\mathcaligr{C}_{r}(W_{n})} I_{y}^{\eta} = 1
\Big|
\sigma\biggr)\nonumber\\[-8pt]\\[-8pt]
& \leq& 1 - \mu_{a,b} ( I_{y}^{\eta} = 1 |
\sigma
) ,\nonumber
\end{eqnarray}
where $\eta$ is an element of $\mathcaligr{C}_{r}(W_{n})$ satisfying
$W_{n}(\eta)=W_{n}$. Lemma \ref{lemm:controlproba} gives a bound for
(\ref{unmoinsprobacond}) which does not depend on the configuration
$\sigma\in\{-,+\}^{\delta B(y,r)}$ nor on the vertex $y$. Hence, it follows
from (\ref{utilisationmarkovian})
\[
\mu_{a,b} ( ^{c}\Gamma_{n} ) \leq( 1 - c_{r} W_{n} )
^{\tau_{n}}
\]
(recall that $1-c_{r}W_{n}$ is positive). Finally, using $\tau_{n}\geq
\tau
R(n)^{d}$ and the classical inequality $\ln(1+X)\leq X$ valid for
$X>-1$, a
bound for the probability of $X_{n}$ to be null is obtained:
\[
\mu_{a,b} (^{c}\Gamma_{n}) \leq\exp( - \tau R(n)^{d} c_{r} W_{n}
) .
\]
Inequality (\ref{stepupperbound}) is proved with $K'=\tau c_{r}$.
\end{pf*}

\section{Distance between local configurations}
\label{section:distancebetween}

\subsection{Motivation}
\label{subsection:motivation}

Throughout Section \ref{section:distancebetween}, $\eta$ represents a
local configuration with radius $r$ having at least one positive
vertex, that is, different from the negative local configuration $\eta
_{-}$. The goal of this section is to describe the model under the
hypothesis (\ref{hypH}), bearing on the potentials $a(n)$ and $b(n)$, and
defined by
\renewcommand{\theequation}{H}
\begin{equation}\label{hypH}
\cases{
n^{d} W_{n}(\eta) = \lambda, &\quad for $\lambda> 0$, \cr
a(n) + 2 \mathcaligr{V} b(n) \leq0, \cr
a(n) + \mathcaligr{V} b(n) \to-\infty, &\quad as $n \to+\infty$.}
\end{equation}
Let us start by giving the reason for (\ref{hypH}). Let $X_{n}(\eta)$
be the
random variable which counts the number of copies of $\eta$ in the
whole graph
$G_{n}$:
\[
X_{n}(\eta) = \sum_{x \in V_{n}} I_{x}^{\eta} .
\]
In \cite{Coupier06-2}, various results have been stated about the
variable $X_{n}(\eta)$. If the product $n^{d}W_{n}(\eta)$ tends to
$0$ (resp., $+\infty$) then the probability $\mu_{a,b}(X_{n}(\eta
)>0)$ tends to $0$ (resp., $1$). In other words, if $W_{n}(\eta)$ is
small compared to $n^{-d}$, then asymptotically, there is no occurrence
of $\eta$ in $G_{n}$. If $W_{n}(\eta)$ is large compared to $n^{-d}$,
then at least one occurrence of $\eta$ can be found in the graph, with
probability tending to $1$. Moreover, provided the hypothesis (\ref
{hypH}) is
satisfied [in particular $W_{n}(\eta)$ is of order $n^{-d}$], the
distribution of $X_{n}(\eta)$ converges weakly to the Poisson
distribution with parameter $\lambda$.

In particular, the number of copies of $\eta$ occurring in the graph
is finite. Let $\eta_{x}$ be one of them (that occurring on the ball
centered at $x$). Lemma \ref{lemm:autour} says vertices around $\eta
_{x}$ are negative with probability tending to $1$. For that purpose,
let us introduce the \textit{ring} $\mathcaligr{R}(x,r,R)$ defined as
the following set of vertices:
\[
\mathcaligr{R}(x,r,R) = \{y \in V_{n} , r < \operatorname{dist}(x,y)
\leq R \} .
\]
\begin{lemm}
\label{lemm:autour}
Let $\eta$ be a local configuration with radius $r$ and let $x\in
V_{n}$. Let us denote by $\mathcaligr{A}_{r}(x,R,+)$ the (interpreted)
event ``there exists at least one positive vertex in $\mathcaligr{R}(x,r,R)$.''

If the potentials satisfy $a(n)+2\mathcaligr{V}b(n)\leq0$ and $a(n)\to
-\infty$. Then,
\[
\lim_{n\to+\infty} \mu_{a,b} \bigl( \mathcaligr
{A}_{r}(x,R,+) |
I_{x}^{\eta}=1 \bigr) = 0.
\]
\end{lemm}
\begin{pf}
Let us introduce the set $\mathcaligr{C}_{R}(\eta,+)$ formed from local
configurations with radius $R$ whose restriction to the ball $B(0,r)$
is equal
to $\eta$ and having at least one positive vertex in the ring
$\mathcaligr{R}(0,r,R)$:
\[
\mathcaligr{C}_{R}(\eta,+) = \{\zeta\in\mathcaligr{C}_{R} ,
\forall
y \in B(0,r) , \zeta(y) = \eta(y) \mbox{ and } k(\zeta
)>k(\eta)
\} .
\]
The intersection of the events $\mathcaligr{A}_{r}(x,R,+)$ and
$I_{x}^{\eta}=1$
forces the ball $B(x,R)$ to contain at least $k(\eta)+1$ positive
vertices. Precisely,
\[
\mathcaligr{A}_{r}(x,R,+) \cap\{ I_{x}^{\eta}=1 \} =
\bigcup_{\zeta\in\mathcaligr{C}_{R}(\eta,+)} \{ I_{x}^{\zeta}=1 \} .
\]
Then, Lemma \ref{lemm:controlproba} allows us to bound the probability
of that
intersection:
%
%
\setcounter{equation}{15}
\renewcommand{\theequation}{\arabic{equation}}
\begin{eqnarray}
\label{smallring}
\mu_{a,b} \bigl( \mathcaligr{A}_{r}(x,R,+) \cap\{ I_{x}^{\eta}=1 \}
\bigr)
&\leq& \sum_{\zeta\in\mathcaligr{C}_{R}(\eta,+)} \mu_{a,b} (
I_{x}^{\zeta}=1
) \nonumber\\[-8pt]\\[-8pt]
&\leq& C_{R} \sum_{\zeta\in\mathcaligr{C}_{R}(\eta,+)}
W_{n}(\zeta).\nonumber
\end{eqnarray}
Let $\zeta$ be an element of $\mathcaligr{C}_{R}(\eta,+)$. By definition,
$k(\zeta)\geq k(\eta)+1$. Moreover, thanks to the convexity of balls, the
perimeter of $\zeta$ is necessary as large as that of $\eta$. In
other words,
$W_{n}(\zeta)\leq\exp(2a(n))W_{n}(\eta)$. Then, we deduce from
(\ref{smallring}) that the conditional probability satisfies
\[
\mu_{a,b} \bigl( \mathcaligr{A}_{r}(x,R,+) |
I_{x}^{\eta}=1
\bigr) \leq c_{r}^{-1} C_{R} |\mathcaligr{C}_{R}(\eta,+)| \exp(
2a(n) )
\]
and tends to $0$ as $n$ tends to infinity.
\end{pf}

This result constitutes the starting point of our study since from now
on, a natural question concerns the distance from $\eta_{x}$ to the
closest positive vertices. Theorems \ref{theo:dist1} and \ref
{theo:dist2} will answer to this question and, more generally, will
describe the distance between local configurations under (\ref{hypH}).

\subsection{The main results}
\label{subsection:main}

The main results of Section \ref{section:distancebetween}, Theorems
\ref{theo:dist1} and \ref{theo:dist2}, concern the two following events.
\begin{defi}
\label{defi:AH}
Let $r\geq1$, $R\geq3r+1$, $0<W<1$ and $x\in V_{n}$. Let us denote by
$\mathcaligr{A}_{r}(x,R,W)$ the event
\[
\exists y \in\mathcaligr{R}(x,2r,R-r) , \exists\zeta\in\mathcaligr
{C}_{r}(W) \qquad\mbox{such that } I_{y}^{\zeta} = 1
\]
and, for all integer $\ell>2r$, by $\mathcaligr{H}_{r}(x,R,\ell,W)$ the event
\begin{eqnarray*}
&&\exists y , y' \in\mathcaligr{R}(x,2r,R-r) , \exists\zeta, \zeta'
\in\mathcaligr{D}_{r}(W)
\\
&&\qquad\mbox{such that }
2r < \operatorname{dist}(y,y') \leq\ell \quad\mbox{and}\quad
I_{y}^{\zeta} =
I_{y'}^{\zeta'} = 1 ,
\end{eqnarray*}
where $\mathcaligr{D}_{r}(W)$ is the set of local configurations with
radius $r$ whose weight is equal to $W$.
\end{defi}

The event $\mathcaligr{A}_{r}(x,R,W)$ means a local configuration with
radius $r$ whose weight is smaller than $W$ can be found somewhere in
the ring $\mathcaligr{R}(x,r,R)$. It generalizes the notation
$\mathcaligr
{A}_{r}(x,R,+)$ of Lemma \ref{lemm:autour}. Indeed, a single positive
vertex can be viewed as a local configuration of weight $W=\exp
(2a(n)-2\mathcaligr{V}b(n))$.

The event $\mathcaligr{H}_{r}(x,R,\ell,W)$ corresponds to the occurrence
in $\mathcaligr{R}(x,r,R)$ of two local configurations of weight $W$ at
distance from each other smaller than $\ell$. The precaution
$\operatorname{dist}(y,y')>2r$ ensures that the occurrences of $\zeta$
and $\zeta'$
do not use the same positive vertices.

The hypothesis (\ref{hypH}) forces the weight of any given local
configurations having at least one positive vertex to tend to $0$ as
$n$ tends to infinity. Hence, the study of the conditional probabilitiy
$\mu_{a,b}(\mathcaligr{A}_{r}(x,R,W)|I_{x}^{\eta}=1)$ becomes relevant
when the weight $W=W_{n}$ is allowed to depend on the size $n$ and to
tend to $0$. In order to avoid trivial situations, it is reasonable to
assume that the set $\mathcaligr{C}_{r}(W_{n})$ is nonempty. Moreover,
remark the events $\mathcaligr{A}_{r}(x,R,W_{n})$ and $\mathcaligr
{A}_{r}(x,R,W_{n}')$ are equal for
\[
W_{n}' = \max_{\zeta\in\mathcaligr{C}_{r}(W_{n})} W_{n}(\zeta) .
\]
As a consequence and without loss of generality, we can assume that,
for each
$n$, the weight $W_{n}$ belongs to $\{ W_{n}(\zeta),
\zeta\in\mathcaligr{C}_{r}\}$. From then on, Lemma \ref{lemm:autour} implies
\[
\mu_{a,b} \bigl( \mathcaligr{A}_{r}(x,R,W_{n})|
I_{x}^{\eta} = 1 \bigr)\to 0,
\]
as $n\to+\infty$, for any fixed radius $R$. So, in order to obtain a
positive limit for the previous quantity, it is needed to take a radius
$R=R(n)$ which tends to $+\infty$. All these remarks are still true
for $\mathcaligr{H}_{r}(x,R,\ell,W)$.

To sum up, the sequences $(R(n))_{n\in\mathbb{N}}$ and $(W_{n})_{n\in
\mathbb{N}}$ will be assumed to satisfy the hypothesis (\ref{hypHprime}):
\renewcommand{\theequation}{H$'$}
\begin{equation}\label{hypHprime}
\cases{
R(n) \to+\infty,\cr
W_{n} \to0, \cr
\forall n , \exists\zeta_{n} \in\mathcaligr{C}_{r} \qquad W_{n} =
W_{n}(\zeta_{n}).}
\end{equation}
Imagine a copy of $\eta$ occurs on $B(x,r)$. Theorem \ref{theo:dist1}
says that the distance to $\eta_{x}$ from the closest positive
vertices is of order
\[
\bigl( \exp\bigl( 2a(n)-2\mathcaligr{V}b(n) \bigr) \bigr) ^{-1/d}
\]
and, more generally, the distance to $\eta_{x}$ from the closest local
configurations of weight $W_{n}$ is of order $W_{n}^{-1/d}$.
\begin{theo}
\label{theo:dist1}
Let us consider a local configuration $\eta\in\mathcaligr{C}_{r}$ and
potentials $a(n)<0$ and $b(n)\geq0$ satisfying (\ref{hypH}). Let
$(R(n))_{n\in\mathbb{N}}$ and $(W_{n})_{n\in\mathbb{N}}$ be two
sequences satisfying (\ref{hypHprime}). Let $x$ be a vertex. Then,
%
%
\setcounter{equation}{16}
\renewcommand{\theequation}{\arabic{equation}}
\begin{eqnarray}
\label{part1}
&&\lim_{n \to+\infty} R(n) W_{n}^{1/d} = 0 \nonumber\\[-8pt]\\[-8pt]
&&\quad \Longrightarrow\quad \lim_{n \to+\infty} \mu_{a,b} \bigl(
\mathcaligr{A}_{r}(x,R(n),W_{n}) | I_{x}^{\eta} = 1 \bigr)
= 0 ;\nonumber
\\
\label{part2}
&&\lim_{n \to+\infty} R(n) W_{n}^{1/d} = +\infty\nonumber\\[-8pt]\\[-8pt]
&&\quad \Longrightarrow\quad \lim_{n \to+\infty} \mu_{a,b} \bigl(
\mathcaligr{A}_{r}(x,R(n),W_{n}) | I_{x}^{\eta} = 1 \bigr)
= 1 .\nonumber
\end{eqnarray}
\end{theo}

First, let us precise that the second part of Theorem \ref{theo:dist1}
concerns only the case where the speed of convergence to $0$ of the
weight $W_{n}$ is slower than that of $W_{n}(\eta)$. Indeed, the
radius $R(n)$ cannot exceed the size $n$ of the graph $G_{n}$ and the
product $n^{d}W_{n}(\eta)$ is constant.

\begin{figure}

\includegraphics{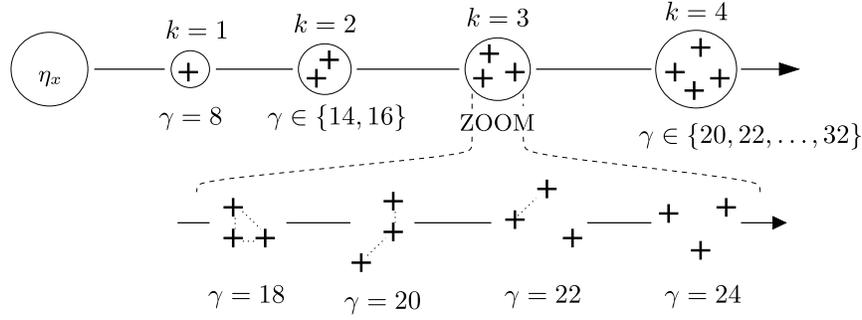}

\caption{Let $\rho=1$, $q=\infty$ and $d=2$; each vertex has
$\mathcaligr{V}=8$ neighbors. Assume potentials satisfy the hypothesis
(\protect\ref{hypH}), $a(n)\to-\infty$ and the ratio $b(n)/a(n)\to0$. A
copy of
$\eta$ occurs on the ball $B(x,r)$. First, the closest local
configurations to $\eta_{x}$ are represented according to their number
of positive vertices. Thus, the situation of those having $k=3$
positive vertices is specified. The closest ones to $\eta_{x}$ are
those having a perimeter $\gamma=18$ and the farthest ones are those
having isolated positive vertices, that is, $\gamma=24$.}
\label{fig:zoom}
\end{figure}

Let us describe Theorem \ref{theo:dist1} when potentials $a(n)$ and
$b(n)$ satisfy the following hypotheses:
\[
a(n) \to-\infty,\qquad b(n) \to+\infty \quad\mbox{and}\quad \frac
{b(n)}{a(n)} \to0 .
\]
Assume a copy of $\eta$ occurs on $B(x,r)$ and consider two local
configurations $\eta_{1}$ and $\eta_{2}$ (whose weights are larger
than that of $\eta$). Then, this is first the number of positive
vertices thus the perimeter which allow to know what are the closest to
$\eta_{x}$, among the copies of $\eta_{1}$ and $\eta_{2}$. Indeed,
if $k(\eta_{1})<k(\eta_{2})$ [$\leq k(\eta)$] then the distance from
the closest copies of $\eta_{1}$ to $\eta_{x}$ is smaller than the
distance from the closest copies of $\eta_{2}$ to $\eta_{x}$,
whatever their perimeters. However, if $\eta_{1}$ and $\eta_{2}$ have
the same number of positive vertices then that having the smallest
perimeter will have its closest copies (to $\eta_{x}$) closer to $\eta
_{x}$ than the closest copies (to $\eta_{x}$) of the other one. Figure
\ref{fig:zoom} proposes a panorama of this situation.

%

%
\begin{figure}[b]

\includegraphics{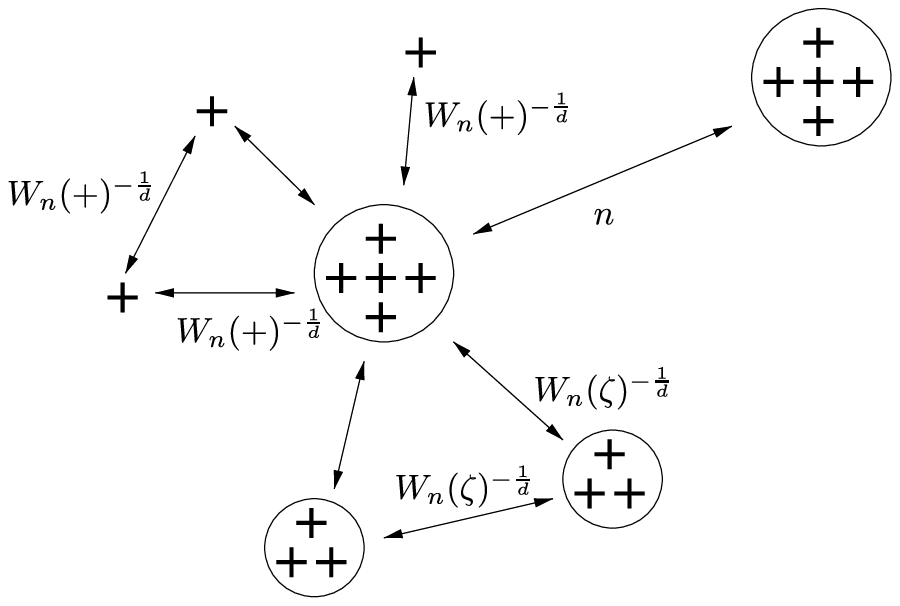}

\caption{Let $\eta$ be the local configuration having $k(\eta)=5$
positive vertices represented above. Assume potentials $a(n)$ and
$b(n)$ satisfy (\protect\ref{hypH}). The distance from a given copy of
$\eta$
occurring in the graph, say $\eta_{x}$, to its closest positive
vertices is of order $W_{n}(+)^{-1/d}$ where $W_{n}(+)$ denotes the
weight of a single positive vertex. The distance between such positive
vertices is of order $W_{n}(+)^{-1/d}$ too. Let $\zeta$ be the local
configuration having $k(\zeta)=3$ positive vertices represented above.
The distance from the closest copies of $\zeta$ to $\eta_{x}$ is of
order $W_{n}(\zeta)^{-1/d}$ and the distance between such copies is
$W_{n}(\zeta)^{-1/d}$ too. Finally, the other copies of $\eta$
occurring in $G_{n}$ are at distance to $\eta_{x}$ of order $n$, since
$n^{d}W_{n}(\eta)=\lambda$.}
\label{fig:dist}
\end{figure}

Now, let us consider two copies of a given local configuration $\zeta$
which are at distance to $\eta_{x}$ of order $W_{n}(\zeta)^{-1/d}$.
Then, the distance between these two copies necessarily tends to
infinity. Otherwise, they would form together a ``super'' local
configuration whose weight would be smaller than $W_{n}(\zeta)$:
thanks to Theorem \ref{theo:dist1}, such an object would not be at
distance to $\eta_{x}$ of order $W_{n}(\zeta)^{-1/d}$. Theorem \ref
{theo:dist2} precises this situation. It says the distance between
these two copies of $\zeta$ tends to infinity as $W_{n}(\zeta
)^{-1/d}$.

Before stating Theorem \ref{theo:dist2}, let us introduce the \textit
{index} of $W_{n}$. The hypothesis (\ref{hypHprime}) ensures the set
$\mathcaligr
{D}_{r}(W_{n})$ is nonempty. However, it is not necessary reduced to
only one element; it may even contain some local configurations not
having the same number of positive vertices. Hence, let us denote by
$k_{n}$ and call the index of $W_{n}$ this maximal number:
\[
k_{n} = \max_{\zeta\in\mathcaligr{D}_{r}(W_{n})} k(\zeta) .
\]
\begin{theo}
\label{theo:dist2}
Let us consider a local configuration $\eta\in\mathcaligr{C}_{r}$ and
potentials $a(n)<0$ and $b(n)\geq0$ satisfying (\ref{hypH}). Let
$(R(n))_{n\in\mathbb{N}}$ and $(W_{n})_{n\in\mathbb{N}}$ be two
sequences satisfying (\ref{hypHprime}). Let $(\ell(n))_{n\in\mathbb
{N}}$ be a
sequence of integers and $x$ be a vertex. Then,
%
%
\setcounter{equation}{18}
\renewcommand{\theequation}{\arabic{equation}}
\begin{eqnarray}
\label{part3}
&&\lim_{n \to+\infty} R(n) W_{n}^{2/d} \max\biggl\{
\ell(n) , \exp\biggl( \frac{2}{d}\mathcaligr{V}b(n)k_{n} \biggr)
\biggr\}= 0 \nonumber\\[-8pt]\\[-8pt]
&&\quad \Longrightarrow\quad \lim_{n \to+\infty} \mu_{a,b} \bigl(
\mathcaligr{H}_{r}(x,R(n),\ell(n),W_{n}) | I_{x}^{\eta} = 1
\bigr) = 0 ;
\nonumber\\
\label{part4}
&&\lim_{n \to+\infty} R(n) \ell(n) W_{n}^{2/d} = +\infty
\nonumber\\[-8pt]\\[-8pt]
&&\quad \Longrightarrow\quad \lim_{n \to+\infty} \mu_{a,b} \bigl(
\mathcaligr{H}_{r}(x,R(n),\ell(n),W_{n}) | I_{x}^{\eta} = 1
\bigr) = 1 .\nonumber
\end{eqnarray}
\end{theo}

Assume a copy of $\eta$ occurs on $B(x,r)$. We know (Theorem
\ref{theo:dist1}) the distance from $\eta_{x}$ to the closest (copies
of) local configurations of weight $W_{n}$ is of order $W_{n}^{-1/d}$.
Let us suppose the radius $R(n)$ is of order $W_{n}^{-1/d}$
(i.e., the product $R_{n}^{d}W_{n}$
tends to a positive constant) and let $\zeta_{y}$ be such a copy [i.e.,
$\zeta\in\mathcaligr{D}_{r}(W_{n})$ occurs on $B(y,r)$ and
$\operatorname{dist}(x,y)$ is
of order $R(n)$]. If $\ell(n)$ is large compared to $W_{n}^{-1/d}$,
then the product $R(n)\ell(n)W_{n}^{2/d}$ tends to infinity and
(\ref{part4}) says a copy of a local configuration of weight $W_{n}$
can be found at distance from $\zeta_{y}$ smaller than $\ell(n)$.
Conversely, if $\ell(n)$ is small compared to $W_{n}^{-1/d}$ then
\[
R(n) W_{n}^{2/d} \max\biggl\{\ell(n) , \exp\biggl( \frac
{2}{d}\mathcaligr{V}b(n)k_{n} \biggr) \biggr\} \to 0 .
\]
Indeed,
\begin{eqnarray*}
W_{n} \exp( 2 \mathcaligr{V} b(n) k_{n} ) & = & W_{n} \exp
( -2 a(n) k_{n} ) \exp\bigl( 2 k_{n} \bigl(a(n) + \mathcaligr
{V} b(n)\bigr) \bigr) \\
& \leq& \exp\bigl( a(n) + \mathcaligr{V} b(n) \bigr)
\end{eqnarray*}
by definition of the index $k_{n}$. So, $W_{n}\exp(2\mathcaligr
{V}b(n)k_{n})$ tends to $0$ as $n\to+\infty$ thanks to the hypothesis
(\ref{hypH}). From then on, (\ref{part3}) says there is no copy of local
configurations of weight $W_{n}$ at distance from $\zeta_{y}$ smaller
than $\ell(n)$. Figure \ref{fig:dist} represents various local
configurations and the distances between each other.

Besides, if the product $R(n)W_{n}^{2/d}$ tends to infinity then (\ref
{part4}) implies that, at distance from $\eta_{x}$ of order $R(n)$,
one can find two copies of local configurations of weight $W_{n}$ so
close to each other that we wish. This is not surprising: Theorem~\ref
{theo:dist1} claims there are local configurations of weight
$W_{n}^{2}$ in the ring $\mathcaligr{R}(x,r,R(n))$ provided
$R(n)W_{n}^{2/d}\to+\infty$.

Finally, let us underline that Theorems \ref{theo:dist1} and \ref
{theo:dist2} would remain unchanged if the radius of local
configurations occurring in $\mathcaligr{R}(x,r,R(n))$ might be different
from that of $\eta$.


Let us end this section by the following remark. Assume that the local
configuration $\eta$ occurs on $B(x,r)$. As long as the ratio $R(n)/n$
tends to $0$, Theorem \ref{theo:dist1} says there is no other copy of
$\eta$ in the ring $\mathcaligr{R}(x,r,R(n))$. Now, $n$ represents the
size of the graph $G_{n}$. So, if other copies (than $\eta_{x}$) of
$\eta$ occur in the graph, they should be at distance of order $n$
from $\eta_{x}$. This is the meaning of Proposition \ref
{prop:dist0(n)}. Recall that $X_{n}(\eta)$ represents the number of
copies of $\eta$ in $G_{n}$. We denote by $\mathcaligr{E}_{n}(\eta)$
the event
\[
\exists C > 0 , ( \forall x \not= y , I_{x}^{\eta} =
I_{y}^{\eta} = 1 ) \quad\Longrightarrow\quad\operatorname{dist}(x,y) \geq C n .
\]
\begin{prop}
\label{prop:dist0(n)}
Let us consider a local configuration $\eta\in\mathcaligr{C}_{r}$ and
potentials $a(n)<0$ and $b(n)\geq0$ satisfying (\ref{hypH}). Then,
\[
\lim_{n \to+\infty} \mu_{a,b} \bigl( \mathcaligr{E}_{n}(\eta
) | X_{n}(\eta) \geq2 \bigr) = 1 .
\]
\end{prop}

Let us recall that the random variable $X_{n}(\eta)$ converges weakly
to a Poisson distribution when the hypothesis (\ref{hypH}) is satisfied (see
\cite{Coupier06-2}). In particular, the fact that no more than one copy
of $\eta$ (i.e., $0$ or $1$) occurs in
the graph has a positive asymptotic probability. So, conditioning by
$X_{n}(\eta)\geq2$, we avoid this uninteresting case.

\subsection{\texorpdfstring{Proofs of Theorems \protect\ref{theo:dist1}, \protect
\ref{theo:dist2} and Proposition \protect\ref{prop:dist0(n)}}{Proofs of
Theorems 5.3, 5.4 and Proposition 5.5}}
\label{subsection:proofs}

The intuition behind Theorem \ref{theo:dist1} is the following. On the
one hand, the Markovian character of the Gibbs measure $\mu_{a,b}$
implies that the events $\mathcaligr{A}_{r}(x,R(n),W_{n})$ and
$\{I_{x}^{\eta}=1\}$ can be considered as asymptotically independent.
So, as $n$ goes to infinity, $\mu_{a,b}(\mathcaligr
{A}_{r}(x,R(n),W_{n}) |
I_{x}^{\eta}=1)$ and $\mu_{a,b}(\mathcaligr{A}_{r}(x,R(n),W_{n}))$ evolve
in the same way. On the other hand, the events
$\mathcaligr{A}_{r}(x,R(n),W_{n})$ and $\mathcaligr{A}(x,R(n),W_{n})$ (see
Section \ref{section:expobounds}) are the same but for a finite number
of vertices, those belonging to $B(x,r)$. Their probabilities have the
same limit. In conclusion, Theorem \ref{theo:expobounds} implies the
conditional probability $\mu_{a,b}(\mathcaligr{A}_{r}(x,R(n),W_{n}) |
I_{x}^{\eta}=1)$ should tend to $0$ or $1$ according to the quantity
$R(n)^{d}W_{n}$. The same remarks hold for Theorem \ref{theo:dist2}.

Thanks to the invariance translation of the graph $G_{n}$, we only
prove Theorems~\ref{theo:dist1} and \ref{theo:dist2} for $x=0$ and we
will denote, respectively, by $\mathcaligr{A}_{r}(R(n),W_{n})$ and
$\mathcaligr{H}_{r}(R(n),\ell(n),W_{n})$ the events
$\mathcaligr{A}_{r}(0,R(n),W_{n})$ and
$\mathcaligr{H}_{r}(0,R(n),\ell(n),W_{n})$. Let us start with the proof of
Theorem \ref{theo:dist1}.
\begin{pf*}{Proof of Theorem \protect\ref{theo:dist1}}
Recall the definition of the event $\mathcaligr{A}_{r}(R(n),W_{n})$:
\[
\exists y \in\mathcaligr{R}\bigl(0,2r,R(n)-r\bigr) , \exists\zeta\in\mathcaligr
{C}_{r}(W_{n}) \qquad I_{y}^{\zeta} = 1 .
\]


Let us start with the proof of (\ref{part1}). The first step is to
move the occurrence of $\zeta$ [the local configuration fulfilling
$\mathcaligr{A}_{r}(R(n),W_{n})$] away from the ball $B(0,r)$ on which
occurs $\eta$. Let $y$ be a vertex belonging to the ring $\mathcaligr
{R}(0,2r,R-r)$. Either $\operatorname{dist}(y,0)\leq2r+1$ and the ball
$B(y,r)$ is
included in $\mathcaligr{R}(0,r,3r+1)$. Either $\operatorname
{dist}(y,0)> 2r+1$ and
$B(y,r)$ is included in $\mathcaligr{R}(0,r+1,R(n))$. In other words,
%
%
\setcounter{equation}{20}
\renewcommand{\theequation}{\arabic{equation}}
\begin{equation}
\label{inclusioneloignement}
\mathcaligr{A}_{r}(R(n),W_{n}) \subset\mathcaligr{A}_{r}(3r+1,W_{n})
\cup
\mathcaligr{A}_{r+1}(R(n),W_{n}) .
\end{equation}
The hypothesis (\ref{hypHprime}) forces the weight $W_{n}$ to be
smaller than
$\exp(2a(n)-2\mathcaligr{V}b(n))$. Hence, the event $\mathcaligr
{A}_{r}(3r+1,W_{n})$ is included in the event $\mathcaligr
{A}_{r}(0,3r+1,+)$ introduced in Section \ref{subsection:motivation}.
Lemma \ref{lemm:autour} implies that the conditional probability
\[
\mu_{a,b} \bigl( \mathcaligr{A}_{r}(3r+1,W_{n}) |
I_{0}^{\eta}=1 \bigr)
\]
tends to $0$. From (\ref{inclusioneloignement}), it remains to prove
the same limit holds for
\[
\mu_{a,b} \bigl( \mathcaligr{A}_{r+1}(R(n),W_{n}) | I_{0}^{\eta} = 1\bigr) .
\]
The strategy consists in introducing increasing events in order to use
the FKG inequality. Since $a(n)+\mathcaligr{V}b(n)$ is negative, the
following events are both increasing:
\[
\bigcup_{\zeta\in\mathcaligr{C}_{r}(W_{n})} \{ I_{y}^{\zeta} = 1 \}
\quad\mbox{and}\quad \bigcup_{\eta' \in\mathcaligr{C}_{r}(W_{n}(\eta
))} \{ I_{0}^{\eta'} = 1 \} .
\]
So does their intersection
\[
S_{y} = \biggl( \bigcup_{\zeta\in\mathcaligr{C}_{r}(W_{n})} \{
I_{y}^{\zeta} = 1 \} \biggr) \cap \biggl( \bigcup_{\eta'
\in\mathcaligr{C}_{r}(W_{n}(\eta))} \{ I_{0}^{\eta'} = 1 \} \biggr) .
\]
As a consequence, we get from the FKG inequality
%
%
\begin{equation}
\label{FKGdist1}\quad
1 - \mu_{a,b} \bigl( \mathcaligr{A}_{r+1}(R(n),W_{n}) \cap
\{ I_{0}^{\eta}=1 \} \bigr) \geq\prod_{y \in
\mathcaligr{R}(0,2r+1,R(n)-r)} \mu_{a,b} ( ^{c}S_{y} ) .
\end{equation}
Let $y\in\mathcaligr{R}(0,2r+1,R(n)-r)$. Let us denote by $V_{y}$ the
union of the two balls $B(y,r)$ and $B(0,r)$, and by $\mathcaligr{C}_{y}$
the set of configurations defined below:
\[
\mathcaligr{C}_{y} = \bigl\{\omega\in\{-,+\}^{V_{y}} ,
W_{n}\bigl(\omega_{B(0,r)}\bigr) \leq W_{n}(\eta) \mbox{ and }
W_{n}\bigl(\omega_{B(y,r)}\bigr) \leq W_{n} \bigr\} .
\]
Then, the probability $\mu_{a,b}(^{c}S_{y})$ becomes
\[
\mu_{a,b} ( ^{c}S_{y} ) = 1 - \sum_{\omega\in\mathcaligr
{C}_{y}} \mu_{a,b} ( I_{V_{y}}^{\omega} = 1 ) ,
\]
where $I_{V_{y}}^{\omega}$ is the random indicator defined on
$\mathcaligr{X}_{n}$ as follows: $I_{V_{y}}^{\omega}(\sigma)$ is $1$ if
the restriction of the configuration $\sigma\in\mathcaligr{X}_{n}$ to
the set $V_{y}$ is $\omega$ and $0$ otherwise. Let $\omega$ be an
element of $\mathcaligr{C}_{y}$. Actually, the inequality (\ref
{controlsupproba}) of Lemma \ref{lemm:controlproba} can be extended to
any subset of vertices $V$ (see Proposition 3.2 of \cite
{Coupier06-2}): whenever $a(n)+2\mathcaligr{V}b(n)$ is negative, there
exists a constant $C_{V}>0$ such that for all $\omega\in\{-,+\}^{V}$,
%
%
\begin{equation}
\label{extension}
\mu_{a,b}(I_{V}^{\omega} = 1) \leq C_{V} W_{n}(\omega) .
\end{equation}
Moreover, $C_{V}$ depends on the set $V$ only through its cardinality.
Applying (\ref{extension}) to $V=V_{y}$, we obtain
\[
\mu_{a,b} ( ^{c}S_{y} ) \geq1 - \sum_{\omega\in
\mathcaligr{C}_{y}} C_{2,r} W_{n}(\omega) ,
\]
where the constant $C_{2,r}$ only depends on $r$. Now, the condition
$\operatorname{dist}(0,y)>2r+1$ forces the connection between the configurations
$\omega_{B(0,r)}$ and $\omega_{B(y,r)}$ to be null. Hence, by the
identity (\ref{lienweights}),
\begin{eqnarray*}
W_{n}(\omega) & = & W_{n}\bigl(\omega_{B(0,r)} \omega_{B(y,r)}\bigr) \\
& = & W_{n}\bigl(\omega_{B(0,r)}\bigr) W_{n}\bigl(\omega_{B(y,r)}\bigr) \\
& \leq& W_{n}(\eta) W_{n} .
\end{eqnarray*}
Therefore, for any vertex $y\in\mathcaligr{R}(0,2r+1,R(n)-r)$,
\[
\mu_{a,b} ( ^{c}S_{y} ) \geq1 - |\mathcaligr{C}_{r}|^{2}
C_{2,r} W_{n}(\eta) W_{n} .
\]
Hence, coupling this latter lower bound with (\ref{FKGdist1}) and
Lemma \ref{lemm:controlproba}, it follows that the conditional probability
\[
\mu_{a,b} \bigl( \mathcaligr{A}_{r+1}(R(n),W_{n}) |
I_{0}^{\eta}=1 \bigr)
\]
is upper bounded by
%
%
\begin{equation}
\label{finalbounddist1}
( c_{r} W_{n}(\eta) )^{-1} \bigl( 1 - \bigl( 1 -
|\mathcaligr{C}_{r}|^{2} C_{2,r} W_{n}(\eta) W_{n}
\bigr)^{(2R(n)+1)^{d}} \bigr) .
\end{equation}
Since the weights $W_{n}$ and $W_{n}(\eta)$ and the product
$R(n)^{d}W_{n}$ tend to $0$, (\ref{finalbounddist1}) is equivalent to
\[
2^{d} c_{r}^{-1} |\mathcaligr{C}_{r}|^{2} C_{2,r} R(n)^{d} W_{n}
\]
and tends to $0$ as $n$ tends to infinity.


Let us turn to the second part of Theorem \ref{theo:dist1}, that is,
statement (\ref{part2}). Some notation introduced in the previous
section will be used here, starting with the set of vertices
\[
T_{n} = \biggl\{ i \bigl(\rho(2r+1)+ 1\bigr) , |i| = 0,\ldots,
\biggl\lfloor
\frac{\rho R(n)}{\rho(2r+1)+ 1} \biggr\rfloor- 1 \biggr\}^{d}
\cap B\bigl(0,R(n)-r\bigr)
\]
(see Figure \ref{fig:reseau}). Let $T_{n}^{\ast}=T_{n}\setminus\{0\}
$ be the same set without the origin. If $\tau_{n}^{\ast}$ denotes
the cardinality of $T_{n}^{\ast}$ then there exists a constant $\tau
^{\ast}>0$ such that $\tau_{n}^{\ast}\geq\tau^{\ast} R(n)^{d}$.
Let us denote by $\mathcaligr{T}_{n}^{\ast}$ the set
\[
\mathcaligr{T}_{n}^{\ast} = \bigcup_{y \in T_{n}^{\ast}} B(y,r)
\]
and by $\Gamma_{n}^{\ast}$ the event
\[
\Gamma_{n}^{\ast} = \bigcup_{y \in T_{n}^{\ast}} \bigcup_{\zeta
\in\mathcaligr{C}_{r}(W_{n})} \{ I_{y}^{\zeta} = 1
\} .
\]
The event $\Gamma_{n}^{\ast}$ implies $\mathcaligr{A}_{r}(R(n),W_{n})$.
Hence, it suffices to prove that the conditional probability $\mu
_{a,b}(^{c}\Gamma_{n}^{\ast} | I_{0}^{\eta}=1)$ tends to $0$. The
end of the proof is very close to that of the upper bound of Theorem
\ref{theo:expobounds}, so it will not be detailed.

The balls forming $\mathcaligr{T}_{n}^{\ast}$ are $\mathcaligr
{V}$-disjoint. So, by the Property \ref{proper:markovian}, the
probability $\mu_{a,b}(X_{n}^{\ast}=0\cap I_{0}^{\eta}=1)$ can be
expressed as
\[
\mathbb{E}_{a,b} \biggl[\mu_{a,b} \bigl( I_{0}^{\eta} = 1
|
\mathcaligr{F}(\delta B(0,r)) \bigr) \prod_{y \in T_{n}^{\ast}} \mu
_{a,b} \biggl( \bigcap_{\zeta\in\mathcaligr{C}_{r}(W_{n})}
I_{y}^{\zeta} = 0 \Big| \mathcaligr{F}(\delta B(y,r)) \biggr)
\biggr] .
\]
Thus, for any vertex $y\in T_{n}^{\ast}$, Lemma \ref
{lemm:controlproba} provides
\[
\mu_{a,b} \biggl( \bigcap_{\zeta\in\mathcaligr{C}_{r}(W_{n})}
I_{y}^{\zeta} = 0 \Big| \mathcaligr{F}(\delta B(y,r)) \biggr)
\leq1 - c_{r} W_{n} .
\]
As a consequence,
\[
\mu_{a,b}( ^{c}\Gamma_{n}^{\ast} \cap I_{0}^{\eta} = 1 ) \leq\mu
_{a,b}( I_{0}^{\eta} = 1 ) ( 1 - c_{r} W_{n} ) ^{\tau
_{n}^{\ast}}
\]
and the conditional probability $\mu_{a,b}(^{c}\Gamma_{n}^{\ast} |
I_{0}^{\eta}=1)$ is now controlled:
\[
\mu_{a,b} ( ^{c}\Gamma_{n}^{\ast} | I_{0}^{\eta} = 1 ) \leq\exp
( - c_{r} \tau^{\ast} R(n)^{d} W_{n} )
\]
which tends to $0$ since by hypothesis $R(n)^{d}W_{n}\to+\infty$.
\end{pf*}

The proof of Theorem \ref{theo:dist2} is based on the same ideas than
that of Theorem \ref{theo:dist1}.
\begin{pf*}{Proof of Theorem \protect\ref{theo:dist2}}
Recall the definition of $\mathcaligr{H}_{r}(R(n),\ell(n),W_{n})$:
\begin{eqnarray*}
&&\exists y , y' \in\mathcaligr{R}\bigl(0,2r,R(n)-r\bigr) , \exists\zeta, \zeta
' \in\mathcaligr{D}_{r}(W_{n}) \\
&&\qquad\mbox{such that }
2r < \operatorname{dist}(y,y') \leq\ell(n) \quad\mbox{and}\quad I_{y}^{\zeta} =
I_{y'}^{\zeta'} = 1 .
\end{eqnarray*}
%

Let us start with the proof of (\ref{part3}). The first step consists
in moving the occurrences of $\eta$, $\zeta$ and $\zeta'$ away from
each other. For that purpose, let us introduce the three following events:
\begin{eqnarray*}
&&\mathcaligr{H}^{1} \dvt \exists y \in V_{n} , \exists\zeta\in
\mathcaligr{D}_{r}(W_{n}) \qquad \operatorname{dist}(0,y)=2r+1
\quad\mbox{and}\quad
I_{y}^{\zeta} = 1 ; \\
&&\mathcaligr{H}^{2} \dvt \exists y , y' \in\mathcaligr{R}\bigl(0,2r,R(n)-r\bigr) ,
\exists\zeta, \zeta' \in\mathcaligr{D}_{r}(W_{n})\\
&&\qquad \mbox{such that }
\operatorname{dist}(y,y') = 2r+1 \quad\mbox{and}\quad I_{y}^{\zeta} =
I_{y'}^{\zeta
'} = 1 ;\\
&&\mathcaligr{H}^{3} \dvt \exists y , y' \in\mathcaligr
{R}\bigl(0,2r+1,R(n)-r\bigr) ,
\exists\zeta, \zeta' \in\mathcaligr{D}_{r}(W_{n})\\
&&\qquad \mbox{such that }
2r+1 < \operatorname{dist}(y,y') \leq\ell(n) \quad\mbox{and}\quad I_{y}^{\zeta}
=I_{y'}^{\zeta'} = 1 ;
\end{eqnarray*}
whose union contains $\mathcaligr{H}_{r}(R(n),\ell(n),W_{n})$. Thanks to
the hypothesis (\ref{hypHprime}), $\mathcaligr{H}^{1}$ is included in
the event
$\mathcaligr{A}_{r}(0,3r+1,+)$ introduced in Section \ref
{subsection:motivation}. So, Lemma \ref{lemm:autour} implies
\[
\lim_{n\to+\infty} \mu_{a,b}(\mathcaligr{H}^{1} | I_{0}^{\eta}=1) =
0 .
\]
Let $y,y'$ be two vertices such that $\operatorname{dist}(y,y')=2r+1$
and $\zeta
,\zeta'$ be two elements of $\mathcaligr{D}_{r}(W_{n})$ occurring on the
balls with radius $r$ centered, respectively, at $y$ and $y'$. The
convexity of balls forces the connection between $\zeta_{y}$ and
$\zeta'_{y'}$ to be smaller than $\frac{1}{2}\mathcaligr{V}k(\zeta)$.
Then, a bound for the weight of the configuration $\zeta_{y}\zeta
'_{y'}$ is deduced:
\begin{eqnarray*}
W_{n}(\zeta_{y}\zeta'_{y'}) & = & W_{n}(\zeta) W_{n}(\zeta') \exp
( 4 b(n) \operatorname{conn}(\zeta_{y},\zeta'_{y'}) ) \\
& \leq& W_{n}^{2} \exp( 2 b(n) \mathcaligr{V}k_{n} ) ,
\end{eqnarray*}
where $k_{n}$ is the index of $W_{n}$. In other words, the event
$\mathcaligr{H}^{2}$ implies the existence in the ring $\mathcaligr
{R}(0,r,R(n))$ of a local configuration with radius $3r+1$ whose weight
is smaller than $W_{n}^{2}\exp(2b(n)\mathcaligr{V}k_{n})$. Now, by hypothesis
\[
R(n)^{d} W_{n}^{2} \exp( 2 b(n) \mathcaligr{V} k_{n} )
\]
tends to $0$. So, it follows from Theorem \ref{theo:dist1}:
\[
\lim_{n\to+\infty} \mu_{a,b}(\mathcaligr{H}^{2} | I_{0}^{\eta}=1) =
0 .
\]
It remains to prove the same limit holds for the event $\mathcaligr
{H}^{3}$. As at the end of the proof of Theorem \ref{theo:dist1}, we
are going to introduce increasing events in order to use the FKG
inequality. For any vertices $y$, $y'$, let us denote by $S_{y,y'}$ the event
\[
\biggl( \bigcup_{\zeta\in\mathcaligr{C}_{r}(W_{n})} \{ I_{y}^{\zeta}
= 1 \} \biggr) \cap \biggl( \bigcup_{\zeta' \in\mathcaligr
{C}_{r}(W_{n})} \{ I_{y'}^{\zeta'} = 1 \} \biggr) \cap
\biggl( \bigcup_{\eta' \in\mathcaligr{C}_{r}(W_{n}(\eta))} \{ I_{0}^{\eta
'} = 1 \} \biggr) .
\]
Since $a(n)+\mathcaligr{V}b(n)$ is negative, $S_{y,y'}$ is increasing.
Thus, the FKG inequality produces the following lower bound:
%
%
\setcounter{equation}{24}
\renewcommand{\theequation}{\arabic{equation}}
\begin{equation}
\label{FKGdist2}
1 - \mu_{a,b} ( \mathcaligr{H}^{3} \cap\{ I_{0}^{\eta
}=1 \} ) \geq\mathop{\mathop{\prod}_{y \in\mathcaligr
{R}(0,2r+1,R(n)-r)}}_{y' \in\mathcaligr{R}(y,2r+1,\ell
(n))} \mu_{a,b} ( ^{c}S_{y,y'} ) .
\end{equation}
Let us pick $y\in\mathcaligr{R}(0,2r+1,R(n)-r)$, $y'\in\mathcaligr
{R}(y,2r+1,\ell(n))$ and denote by $V_{y,y'}$ the union of the three
balls $B(y,r)$, $B(y',r)$ and $B(0,r)$. Let $\mathcaligr{C}_{y,y'}$ be
the set of configurations defined by
\[
\mathcaligr{C}_{y,y'} = \left\{\omega\in\{-,+\}^{V_{y,y'}} ,
\matrix{
W_{n}\bigl(\omega_{B(0,r)}\bigr) \leq W_{n}(\eta) \cr
W_{n}\bigl(\omega_{B(y,r)}\bigr) \leq W_{n} \cr
W_{n}\bigl(\omega_{B(y',r)}\bigr) \leq W_{n}}
\right\}.
\]
Then, using inequality (\ref{controlsupproba}) of Lemma \ref
{lemm:controlproba} [or rather its extension (\ref{extension})], we
bound the probability of $S_{y,y'}$. There exists a constant
$C_{3,r}>0$ such that
\begin{eqnarray*}
\mu_{a,b} ( S_{y,y'} ) & = & \sum_{\omega\in\mathcaligr
{C}_{y,y'}} \mu_{a,b} ( I_{V_{y,y'}}^{\omega} = 1 ) \\
& \leq& C_{3,r} \sum_{\omega\in\mathcaligr{C}_{y,y'}} W_{n}(\omega) .
\end{eqnarray*}
Now, vertices $0$, $y$ and $y'$ are sufficiently far apart so that the
weight of $\omega$ might write as the product of the weights of
$\omega_{B(0,r)}$, $\omega_{B(y,r)}$ and $\omega_{B(y',r)}$. Then,
\[
\mu_{a,b} ( S_{y,y'} ) \leq C_{3,r} |\mathcaligr
{C}_{r}|^{3} W_{n}(\eta) W_{n}^{2}
\]
since $\omega\in\mathcaligr{C}_{y,y'}$ and the cardinality of
$\mathcaligr
{C}_{y,y'}$ is bounded by $|\mathcaligr{C}_{r}|^{3}$. Finally, coupling
this latter inequality with (\ref{FKGdist2}) and Lemma \ref
{lemm:controlproba}, we bound the conditional probability $\mu
_{a,b}(\mathcaligr{H}^{3} | I_{0}^{\eta}=1)$ by
\[
( c_{r} W_{n}(\eta) )^{-1} \bigl( 1 - \bigl( 1 -
|\mathcaligr{C}_{r}|^{3} C_{3,r} W_{n}(\eta) W_{n}^{2}
\bigr)^{(2R(n)+1)^{d}(2\ell(n)+1)^{d}} \bigr) .
\]
Since $W_{n}$, $W_{n}(\eta)$ and the product $\ell
(n)^{d}R(n)^{d}W_{n}^{2}$ tend to $0$, the above quantity is equivalent to
\[
2^{2d} c_{r}^{-1} |\mathcaligr{C}_{r}|^{3} C_{3,r} \ell(n)^{d} R(n)^{d}
W_{n}^{2}
\]
and tends to $0$ as $n$ tends to infinity. So does $\mu
_{a,b}(\mathcaligr
{H}^{3} | I_{0}^{\eta}=1)$.


Let us deal with the second part of Theorem \ref{theo:dist2}, that is,
statement (\ref{part4}).

If the sequence $(\ell(n))_{n\in\mathbb{N}}$ were bounded, say by a
constant $\ell(\infty)$, the event $\mathcaligr{H}_{r}(R(n),\ell
(n),W_{n})$ would correspond to the existence in the ring $\mathcaligr
{R}(0,r,R(n))$ of a local configuration with radius $\ell(\infty)$,
say $\zeta$, and whose weight would be larger than $W_{n}^{2}$ [by
(\ref{lienweights})]. Hence, the limit
\[
R(n)^{d} W_{n}(\zeta) \geq\frac{1}{\ell(\infty)^{d}} R(n)^{d} \ell
(n)^{d} W_{n}^{2} \to +\infty
\]
would imply [statement (\ref{part2}) of Theorem \ref{theo:dist1}] a
limit equal to $1$ for the conditional probability $\mu
_{a,b}(\mathcaligr
{H}_{r}(R(n),\ell(n),W_{n}) | I_{0}^{\eta}=1)$. So, from now on, we
will assume that $\ell(n)\to+\infty$.

The second part of Theorem \ref{theo:dist2} needs one more time the
set of vertices
\[
T_{n} = \biggl\{ i \bigl(\rho(2r+1)+ 1\bigr) , |i| = 0,\ldots,
\biggl\lfloor
\frac{\rho R(n)}{\rho(2r+1)+ 1} \biggr\rfloor- 1 \biggr\}^{d}
\cap B\bigl(0,R(n)-r\bigr)
\]
(see Figure \ref{fig:reseau}). Let $\phi_{n}=i_{n}(\rho(2r+1)+ 1)$, where
\[
i_{n} = \biggl\lfloor\frac{\rho({4}/{3}\ell(n)+1) + 1}{\rho
(2r+1)+ 1} \biggr\rfloor+ 1 ,
\]
and let us denote by $L_{n}$ the sublattice of $T_{n}$ defined by
\[
L_{n} = \bigl\{ ( \phi_{n} \mathbb{Z} ) ^{d}
\cap B \bigl( 0 , R(n)-\tfrac{2}{3}\ell(n) \bigr) \bigr\}
\setminus\{ 0 \} .
\]
Any two balls $B(y,\frac{2}{3}\ell(n))$ and $B(y',\frac{2}{3}\ell
(n))$, where $y,y'\in L_{n}$, are $\mathcaligr{V}$-disjoint since
\[
\operatorname{dist}(y,y') \geq\frac{\phi_{n}}{\rho} > \frac{4}{3} \ell
(n) + 1 .
\]
Obviously, they are also included in the large ball $B(0,R(n))$. Let us
pick a vertex $y$ of $L_{n}$. We denote by $\Lambda_{1}(y)$ and
$\Lambda_{2}(y)$ the following subsets of $T_{n}$:
\[
\Lambda_{1}(y) = \bigl\{ z \in T_{n} , \operatorname{dist}(z,y) \leq\tfrac
{1}{3}\ell(n) \bigr\}
\]
and
\[
\Lambda_{2}(y) = \bigl\{ z \in T_{n} , \tfrac{1}{3}\ell(n) <
\operatorname{dist}(z,y) \leq\tfrac{2}{3}\ell(n) \bigr\}.
\]
Hence, two vertices $z\in\Lambda_{1}(y)$ and $z'\in\Lambda_{2}(y)$
satisfy $2r<\operatorname{dist}(z,z')\leq\ell(n)$. As a consequence,
the event
\[
\Gamma= \bigcup_{y\in L_{n}} \mathop{\mathop{\bigcup}_{z\in\Lambda
_{1}(y)}}_{z'\in\Lambda_{2}(y)} \bigcup_{\zeta
,\zeta'\in\mathcaligr{D}_{r}(W_{n})} \{ I_{z}^{\zeta} =
I_{z'}^{\zeta'} = 1 \}
\]
implies $\mathcaligr{H}_{r}(R(n),\ell(n),W_{n})$. So, it suffices to
prove the conditional probability $\mu_{a,b}(^{c}\Gamma| I_{0}^{\eta
}=1)$ tends to $0$ as $n\to+\infty$.

The vertices of $\Lambda_{1}(y)$ are the elements of $((\rho
(2r+1)+1)\mathbb{Z})^{d}$ whose distance (with respect to \mbox{$\|\cdot\|
_{\infty}$}) to $y$ is smaller than $\frac{\rho}{3}\ell(n)$. Their
number is of order $\ell(n)^{d}$. More generally, there exists a
constant $\alpha$ which depends only on the parameters $d$, $\rho$,
$q$ and $r$, such that, for any vertex $y$,
\[
|\Lambda_{i}(y)| \geq\alpha\ell(n)^{d} ,\qquad i = 1,2\quad \mbox{and}
\quad |L_{n}| \geq\alpha\frac{R(n)^{d}}{\ell(n)^{d}} .
\]
Moreover, the Markovian character of the measure $\mu_{a,b}$ (see
Property \ref{proper:markovian}) allows us to express the probability
$\mu_{a,b}(^{c}\Gamma\cap I_{0}^{\eta}=1)$ as the following expectation:
\begin{eqnarray*}
&&\mathbb{E}_{a,b} \biggl[\mu_{a,b} \bigl(
I_{0}^{\eta} = 1
|\mathcaligr{F}(\delta B(0,r)) \bigr) \nonumber\\[-8pt]\\[-8pt]
&&\qquad{} \times\prod_{y\in L_{n}} \mathop{\mathop{\prod}_{z\in
\Lambda_{1}(y)}}_{z'\in\Lambda_{2}(y)} \mu_{a,b}
\biggl( \bigcap_{\zeta,\zeta'\in\mathcaligr{D}_{r}(W_{n})}
I_{z}^{\zeta} I_{z'}^{\zeta'} = 0 \Big| \mathcaligr{F}(\delta
V_{z,z'}) \biggr) \biggr],\nonumber
\end{eqnarray*}
where $V_{z,z'}$ is the union of the two balls $B(z,r)$ and $B(z',r)$.
Let $y\in L_{n}$ and $z\in\Lambda_{1}(y)$, $z'\in\Lambda_{2}(y)$.
Let $\sigma\in\{-,+\}^{V_{z,z'}}$. Since the balls $B(z,r)$ and
$B(z',r)$ are $\mathcaligr{V}$-disjoint,
\begin{eqnarray*}
&&\mu_{a,b} \biggl( \bigcap_{\zeta,\zeta'\in
\mathcaligr{D}_{r}(W_{n})} I_{z}^{\zeta} I_{z'}^{\zeta'} = 0
\Big| \sigma\biggr) \nonumber\\
&&\qquad = 1 - \sum_{\zeta,\zeta'\in\mathcaligr{D}_{r}(W_{n})} \mu_{a,b}
\bigl( I_{z}^{\zeta} = 0 | \sigma_{B(z,r)} \bigr) \mu_{a,b} \bigl( I_{z'}^{\zeta'}
= 0 | \sigma_{B(z',r)} \bigr) \\
&&\qquad \leq1 - c_{r}^{2} W_{n}^{2} .
\end{eqnarray*}
In conclusion,
\begin{eqnarray*}
\mu_{a,b}(^{c}\Gamma| I_{0}^{\eta}=1) & \leq& ( 1 - c_{r}^{2}
W_{n}^{2} ) ^{\alpha^{3} R(n)^{d} \ell(n)^{d}} \\
& \leq& \exp( - c_{r}^{2} \alpha^{3} R(n)^{d} \ell(n)^{d}
W_{n}^{2} ),
\end{eqnarray*}
which tends to $0$ whenever $R(n)^{d}\ell(n)^{d}W_{n}^{2}$ tends to $0$.
\end{pf*}

Let us finish by the proof of Proposition \ref{prop:dist0(n)}. It is
obtained using the results of the proofs of Lemma \ref{lemm:autour}
and Theorem \ref{theo:dist1}.
\begin{pf*}{Proof of Proposition \protect\ref{prop:dist0(n)}}
We are going to prove that the quantity $\mu_{a,b}(^{c}\mathcaligr
{E}_{n}(\eta) | X_{n}(\eta)\geq2)$ tends to $0$ as $n\to+\infty$,
where the event $^{c}\mathcaligr{E}_{n}(\eta)$ means
\[
\forall C > 0 , \exists x \not= y \qquad I_{x}^{\eta} I_{y}^{\eta}
= 1 \quad\mbox{and}\quad \operatorname{dist}(x,y) < C n .
\]
Thanks to the invariance translation of the graph $G_{n}$ and for any
constant $C>0$,
\begin{eqnarray*}
&& \mu_{a,b} \bigl( ^{c}\mathcaligr{E}_{n}(\eta)
| X_{n}(\eta) \geq2 \bigr) \nonumber\\
&&\qquad \leq n^{d} \mu_{a,b} \bigl( \exists x \not= 0 ,
I_{x}^{\eta} I_{0}^{\eta} = 1 , \operatorname{dist}(x,0) \leq\lfloor C
n \rfloor
| X_{n}(\eta) \geq2 \bigr) .
\end{eqnarray*}
The condition $X_{n}(\eta)\geq2$ can be removed. Indeed, the random
variable $X_{n}(\eta)$ is asymptotically Poissonian, under the
hypothesis (\ref{hypH}) (see \cite{Coupier06-2}). So, there exist a constant
$\varepsilon>1$ and an integer $n(\varepsilon)$ such that the
probability $\mu_{a,b}(X_{n}(\eta)\geq2)$ is larger than
$\varepsilon^{-1}$ whenever $n\geq n(\varepsilon)$. Hence, for such
value of $n$,
\begin{eqnarray*}
&&\mu_{a,b} \bigl( \exists x \not= 0 ,
I_{x}^{\eta} I_{0}^{\eta} = 1 , \operatorname{dist}(x,0) \leq\lfloor C
n \rfloor
| X_{n}(\eta) \geq2 \bigr)\\
&&\qquad \leq\varepsilon\mu_{a,b} \bigl( \exists x \not= 0 ,
I_{x}^{\eta} I_{0}^{\eta} = 1 , \operatorname{dist}(x,0) \leq\lfloor C
n \rfloor
\bigr) .
\end{eqnarray*}
In a first time, assume the centers $0$ and $x$ of balls with radius
$r$ on which the local configuration $\eta$ occurs satisfy
$\operatorname{dist}(x,0)\leq2r+1$. A configuration of $\{-,+\}
^{B(0,R)}$, $R=3r+1$,
fulfilling this event necessarily belongs to
\[
\mathcaligr{C}_{R}(\eta,+) = \{\zeta\in\mathcaligr{C}_{R} ,
\forall
y \in B(0,r) , \zeta(y) = \eta(y) \mbox{ and } k(\zeta
)>k(\eta)
\} .
\]
Now, using the inequalities of the proof of Lemma \ref{lemm:autour},
we get
\begin{eqnarray*}
&&n^{d} \mu_{a,b} \bigl( \exists x \not= 0 , I_{x}^{\eta
} I_{0}^{\eta} = 1 , \operatorname{dist}(x,0) \leq2r+1 \bigr) \\
&&\qquad \leq C_{R} |\mathcaligr{C}_{R}(\eta,+)| n^{d} W_{n}(\eta) \exp
( 2a(n) ),
\end{eqnarray*}
which tends to $0$ since $C_{R}$ and $|\mathcaligr{C}_{R}(\eta,+)|$ are
some constants, $n^{d}W_{n}(\eta)=\lambda$ and the magnetic field
$a(n)\to-\infty$. So, the case where $\operatorname{dist}(x,0)\leq2r+1$ is
negligible and
\begin{eqnarray*}
&&\limsup_{n\to+\infty} \mu_{a,b} \bigl(
^{c}\mathcaligr{E}_{n}(\eta) | X_{n}(\eta) \geq2 \bigr)\\
&&\qquad \leq\varepsilon\limsup_{n\to+\infty} n^{d} \mu_{a,b} \bigl(
\exists x \not= 0 , I_{x}^{\eta} I_{0}^{\eta} = 1 , 2r+1 <
\operatorname{dist}(x,0) \leq\lfloor C n \rfloor\bigr)
\end{eqnarray*}
for any $C>0$. Remark the event
\[
\exists x \not= 0 \qquad I_{x}^{\eta} I_{0}^{\eta} = 1 ,\qquad 2r+1 <
\operatorname{dist}(x,0) \leq\lfloor C n \rfloor
\]
is included in $\mathcaligr{A}_{r+1}(0,\lfloor Cn\rfloor,W_{n}(\eta
))\cap\{I_{0}^{\eta}=1\}$. The proof of Theorem \ref{theo:dist1}
gave us an upper bound for the probability
\[
\mu_{a,b} \bigl( \mathcaligr{A}_{r+1} ( 0, \lfloor Cn \rfloor,
W_{n}(\eta) ) \cap\{ I_{0}^{\eta} = 1 \} \bigr),
\]
which behaves as $\lfloor Cn\rfloor^{d} W_{n}(\eta)^{2}$. Since the
product $n^{d} W_{n}(\eta)$ is constant,
\[
n^{d} \mu_{a,b} \bigl( \exists x \not= 0 , I_{x}^{\eta}
I_{0}^{\eta} = 1 , 2r+1 < \operatorname{dist}(x,0) \leq\lfloor C n
\rfloor\bigr)
\]
is upperbounded by a constant term multiplied by $C^{d}$. Take
$C\searrow0$ and the desired result follows:
\[
\limsup_{n\to+\infty} \mu_{a,b} \bigl( ^{c}\mathcaligr
{E}_{n}(\eta) | X_{n}(\eta) \geq2 \bigr) = 0 .
\]
\upqed\end{pf*}

\section{Ubiquity of local configurations}
\label{section:ubiquity}

This section proposes an application of inequalities stated in the
proof of Theorem \ref{theo:expobounds}. Its goal is to give a
criterion ensuring that a given local configuration occurs everywhere
in the graph. This criterion is presented in Proposition \ref
{prop:ubiquity} through the use of an appropriate Markovian measure
built from the Gibbs measure $\mu_{a,b}$ of (\ref{def:mesgibbs}).

A simple way to cover the set of vertices $V_{n}$ by balls consists in
using the $L_{\infty}$ norm. In this case, replacing the radius $r$
with $\rho r$, we can assume that $\rho=1$. So, in this section, the
neighborhood structure of each vertex is
\[
\mathcaligr{V}(x) = \{ y \in V_{n} , \|y-x\|_{\infty} = 1 \} .
\]
Hence, the graph distance $\mathrm{dist}$ and the $L_{\infty}$ norm
define the
same sets; the ball $B(x,r)$ is equal to the hypercube with center $x$
and radius $r$.

For all integer $n$, let us denote by $v(n)$ the largest integer $m$
such that $2^{m}$ divides $n$ and by $R(n)$ the integer satisfying the relation
\[
\frac{n}{2^{v(n)}} = 2 R(n) + 1 .
\]
A function $g\dvtx\mathbb{N}^{\ast}\to\mathbb{N}^{\ast}$ satisfying
for all integer $n$ the inequalities $v(g(n+1))\geq v(g(n))$ and
$R(g(n+1))\geq R(g(n))$, and tending to infinity as $n\to+\infty$
will be said \textit{adequate}. In particular, an adequate function is
nondecreasing. The functions recursively defined by
\[
\cases{
g(n+1) = 2^{v(g(n))+1} \bigl( 2\bigl(R(g(n))+1\bigr) + 1 \bigr) ,&\quad $\forall n
\geq1$, \cr
g(1) \in\mathbb{N}^{\ast},}
\]
provide examples of adequate functions, since $v(g(n+1))$ and
$R(g(n+1))$ are, respectively, equal to $v(g(n))+1$ and $R(g(n))+1$. For
instance, if $g(1)=1$ then $g(2)=2^{1}\times3=6$, $g(3)=2^{2}\times
5=20$, $g(4)=2^{3}\times7=56,\ldots.$

In conclusion, replacing $n$ with $g(n)$ where $g$ is an adequate
function, we will assume that there exists a nondecreasing, integer
valued sequence $(R(n))_{n\geq1}$ such that the sequence
\[
\biggl( \frac{n}{2R(n)+1} \biggr) _{n \geq1}
\]
is nondecreasing and integer valued too.

Now, let us denote by $\tilde{V}_{n}$ the following subset of $V_{n}$:
\[
\tilde{V}_{n} = \biggl\{ i \bigl(2R(n)+1\bigr) , i = 0,\ldots, \frac
{n}{2R(n)+1} - 1 \biggr\}^{d} .
\]
The set of balls $\{B(x,R(n)), x\in\tilde{V}_{n}\}$ is a partition of
$V_{n}$. The edge set $\tilde{E}_{n}$ is specified by defining the
neighborhood $\tilde{\mathcaligr{V}}(x)$ of each vertex $x\in\tilde{V}_{n}$:
%
%
\setcounter{equation}{25}
\renewcommand{\theequation}{\arabic{equation}}
\begin{equation}
\label{voisins_tilde}
\tilde{\mathcaligr{V}}(x) = \{ y \in\tilde{V}_{n} , \|y-x\|_{\infty
} = 2R(n)+1 \} .
\end{equation}
By analogy with the previous sections, we denote by $\delta\tilde{V}$
the neighborhood of $\tilde{V}\subset\tilde{V}_{n}$ corresponding to
(\ref{voisins_tilde}). Hence, an undirected graph structure $\tilde
{G}_{n}=(\tilde{V}_{n},\tilde{E}_{n})$ with periodic boundary
conditions is defined. The \textit{size} of $\tilde{G}_{n}$ is the
ratio $n$ divided by $2R(n)+1$. Furthermore, remark the graph $\tilde
{G}_{n}$ retains the translation invariance property of $G_{n}$.

Let $\eta$ be a local configuration with radius $r$. We associate to
$\eta$ a function $f$ from $\mathcaligr{X}_{n}$ into $\tilde{\mathcaligr
{X}}_{n}=\{-,+\}^{\tilde{V}_{n}}$ defined for $\sigma\in\mathcaligr
{X}_{n}$ and $x\in\tilde{V}_{n}$ by
\[
f(\sigma)(x) = + \quad\Longleftrightarrow \quad\exists y \in B(x,R(n))
\qquad
I_{y}^{\eta}(\sigma) = 1 .
\]
In other words, the vertex $x\in\tilde{V}_{n}$ is positive for the
configuration $f(\sigma)\in\tilde{\mathcaligr{X}}_{n}$ if and only if
a copy of $\eta$ occurs in the ball $B(x,R(n)+r)$ (under $\sigma$).
Remark the function $f$ is onto $\tilde{\mathcaligr{X}}_{n}$. In
particular, for $\tilde{V}\subset\tilde{V}_{n}$, the values taken by
a configuration $\sigma\in\mathcaligr{X}_{n}$ on the set
\[
\Sigma= \bigcup_{y \in\tilde{V}} B\bigl(y,R(n)+r\bigr)
\]
specify completely the configuration $f(\sigma)$ on $\tilde{V}$.
Besides, $f$ is not an injective function; for $\zeta\in\{-,+\}
^{\tilde{V}}$, $f^{-1}(\zeta)$ represents a subset of $\{-,+\}
^{\Sigma}$.

Let $\mathcaligr{F}(\tilde{V})$ be the $\sigma$-algebra generated by
the configurations of $\{-,+\}^{\tilde{V}}$. Then, $f^{-1}(\mathcaligr
{F}(\tilde{V}))$ is still a $\sigma$-algebra, generated by the sets
$f^{-1}(\zeta)$, $\zeta\in\{-,+\}^{\tilde{V}}$, and is coarser than
$\mathcaligr{F}(\Sigma)$:
\[
f^{-1}(\mathcaligr{F}(\tilde{V})) \varsubsetneq\mathcaligr{F}(\Sigma) .
\]
Thus, let us endow the set of configurations $\tilde{\mathcaligr
{X}}_{n}$ with the measure $\tilde{\mu}_{a,b}$ defined by
%
%
\begin{equation}
\label{mu_tilde}
\forall\zeta\in\tilde{\mathcaligr{X}}_{n} \qquad \tilde{\mu
}_{a,b}(\zeta) = \mu_{a,b}(f^{-1}(\zeta)) .
\end{equation}
Expectations relative to $\tilde{\mu}_{a,b}$ will be denoted by
$\tilde{\mathbb{E}}_{a,b}$.

Property \ref{proper:markovian_tilde} links the random variable
$\tilde{\mu}_{a,b}(\cdot|\mathcaligr{F}(\tilde{V}))$ to $\mu
_{a,b}(\cdot|f^{-1}(\mathcaligr{F}(\tilde{V})))$ and states the
Markovian character of $\tilde{\mu}_{a,b}$. The identity (\ref
{mu_tilde_mu}) is completely based on the definition of the measure
$\tilde{\mu}_{a,b}$. It holds whatever the function $f$. Relation
(\ref{markovian_tilde}) derives from the Markovian character of the
Gibbs measure $\mu_{a,b}$ and the use of the $L_{\infty}$ norm in the
construction of the graph $\tilde{G}_{n}$. Property \ref
{proper:markovian_tilde} will be proved at the end of the section.
\begin{proper}
\label{proper:markovian_tilde}
Let $\tilde{U},\tilde{V}$ be two subsets of $\tilde{V}_{n}$. Then,
for all event $A\in\mathcaligr{F}(\tilde{U})$,
%
%
\begin{equation}
\label{mu_tilde_mu}
\tilde{\mu}_{a,b} ( A | \mathcaligr{F}(\tilde
{V}) ) \circ f = \mu_{a,b} ( f^{-1}(A)
| f^{-1}(\mathcaligr{F}(\tilde{V})) ) ,
\end{equation}
where $\circ$ denotes the composition relation. Moreover, if $\tilde
{U}\cap\tilde{V}=\varnothing$ and $\delta\tilde{U}\subset\tilde
{V}$ then, for all event $A\in\mathcaligr{F}(\tilde{U})$,
%
%
\begin{equation}
\label{markovian_tilde}
\tilde{\mu}_{a,b} ( A | \mathcaligr{F}(\tilde
{V}) ) = \tilde{\mu}_{a,b} ( A |
\mathcaligr{F}(\delta\tilde{U}) ) .
\end{equation}
\end{proper}

The rest of this section is devoted to the study of the asymptotic
behavior of the probability measure $\tilde{\mu}_{a,b}$. As $n$ tends
to infinity, the two sequences
\[
(R(n))_{n\geq1} \quad\mbox{and}\quad \biggl( \frac{n}{2R(n)+1}
\biggr) _{n \geq1}
\]
cannot be simultaneously bounded. Then, three alternatives are
conceivable; either the graph $G_{n}$ is divided into a small number of
large balls, or into a large number of small balls or a large number of
large balls.

Proposition \ref{prop:ubiquity} gives sufficient conditions describing
the asymptotic behavior of $\tilde{\mu}_{a,b}$. Relations (\ref
{zeta-1}) and (\ref{zeta-0}) are a rewriting of results already known;
the second one means at least one copy of the local configuration $\eta
$ can be found somewhere in the graph, with probability tending to $1$
(for $\mu_{a,b}$), whenever $W_{n}(\eta)$ is larger than $n^{-d}$.
Besides, recall that any given ball $B(x,R(n)+r)$ contains a copy of
$\eta$ whenever $W_{n}(\eta)$ is larger than $R(n)^{-d}$ (Theorem
\ref{theo:expobounds}). But, so as to every ball $B(x,R(n)+r)$, $x\in
\tilde{V}_{n}$, contains a copy of $\eta$ (this is that we call
\textit{ubiquity} of the local configuration $\eta$) a stronger
condition is given by Proposition \ref{prop:ubiquity}; $W_{n}(\eta)$
must be larger than $R(n)^{-d} \ln(n/R(n))$.

For that purpose, let $\zeta_{-}$ and $\zeta_{+}$ be the two
configurations of $\tilde{\mathcaligr{X}}_{n}$ whose vertices are all
negative and all positive.
\begin{prop}
\label{prop:ubiquity}
Let us consider a local configuration $\eta\in\mathcaligr{C}_{r}$ and
potentials $a(n)<0$ and $b(n)\geq0$ such that $a(n)+2\mathcaligr
{V}b(n)\leq0$.
%
%
\begin{eqnarray}
\label{zeta-1}
&&\hspace*{22pt}\mbox{If } \lim_{n \to+\infty} n^{d} W_{n}(\eta) = 0
\qquad\mbox{then } \lim_{n \to+\infty} \tilde{\mu}_{a,b} ( \zeta_{-}
) = 1 .
\\
\label{zeta-0}
&&\hspace*{22pt}\mbox{If } \lim_{n \to+\infty} n^{d} W_{n}(\eta) = +\infty
\qquad\mbox{then } \lim_{n \to+\infty} \tilde{\mu}_{a,b} (
\zeta_{-} ) = 0 .
\\
\label{zeta+1}
&&\hspace*{22pt}\mbox{If } \lim_{n \to+\infty} R(n)^{d} \ln\biggl( \frac
{n}{R(n)} \biggr)^{-1} W_{n}(\eta) = +\infty \qquad\mbox{then }
\lim_{n \to+\infty} \tilde{\mu}_{a,b} ( \zeta_{+} )
= 1 .
\end{eqnarray}
\end{prop}

Relations (\ref{zeta-1}) and (\ref{zeta+1}), respectively, say that
$\tilde{\mu}_{a,b}$ converges weakly to the Dirac measures associated
to the configurations $\zeta_{-}$ and $\zeta_{+}$.

The probability (for $\tilde{\mu}_{a,b}$) for a given vertex
$x\in\tilde{V}_{n}$ to be positive is equal to the probability (for
$\mu_{a,b}$) that the local configuration $\eta$ occurs somewhere in
the ball $B(x,R(n)+r)$; see relation (\ref{mutildemu}) below. Now, this
quantity has been studied and bounded in Section
\ref{section:expobounds}. The proof of Proposition \ref{prop:ubiquity}
immediately derives from this remark.
\begin{pf*}{Proof of Proposition \protect\ref{prop:ubiquity}}
For $x\in\tilde{V}_{n}$, let $\tilde{I}_{x}^{+}$ be the indicator
function defined on $\tilde{\mathcaligr{X}}_{n}$ as follows: $\tilde
{I}_{x}^{+}(\zeta)$ is $1$ if the vertex $x$ is positive under $\zeta
\in\tilde{\mathcaligr{X}}_{n}$, that is, $\zeta(x)=1$, and $0$
otherwise. Then,
\begin{eqnarray*}
\tilde{\mu}_{a,b} ( \exists x \in\tilde{V}_{n} , \tilde
{I}_{x}^{+} = 1 ) & = & \mu_{a,b} ( \exists x \in V_{n} ,
I_{x}^{\eta} = 1 ) \\
& = & \mu_{a,b} \bigl( X_{n}(\eta) > 0 \bigr) ,
\end{eqnarray*}
where $X_{n}(\eta)$ represents the number of copies of $\eta$
occurring in $G_{n}$. In \cite{Coupier06-2}, it has been proved that
$\mu_{a,b}(X_{n}(\eta)>0)$ tends to $0$ (resp., $1$) whenever
$n^{d}W_{n}(\eta)$ tends to $0$ (resp., $+\infty$). Relations (\ref
{zeta-1}) and (\ref{zeta-0}) follow.

In order to obtain (\ref{zeta+1}), we are going to prove that the
probability of the opposit event $\exists x\in\tilde{V}_{n}, \tilde
{I}_{x}^{+}=0$ tends to $0$. Due to the translation invariance of
$\tilde{G}_{n}$, the random indicators $\tilde{I}_{x}^{+}$ have the
same distribution. So,
\[
\tilde{\mu}_{a,b} ( \exists x \in\tilde{V}_{n} , \tilde
{I}_{x}^{+} = 0 ) \leq\biggl(\frac{n}{R(n)} \biggr)^{d} \tilde
{\mu}_{a,b} ( \tilde{I}_{0}^{+} = 0 ) .
\]
Furthermore, following the inequalities of the proof of the upper bound
of (\ref{expobounds}), we bound the probability of the event $\tilde
{I}_{0}^{+}=0$,
%
%
\begin{eqnarray}
\label{mutildemu}
\tilde{\mu}_{a,b} ( \tilde{I}_{0}^{+} = 0 ) & = & \mu
_{a,b} ( \tilde{I}_{0}^{+} \circ f = 0 ) \nonumber\\
& = & \mu_{a,b} \biggl( \sum_{y \in B(0,R(n))} I_{y}^{\eta} = 0
\biggr) \\
& \leq& \exp( - \tau c_{r} R(n)^{d} W_{n}(\eta) ) ,
\nonumber
\end{eqnarray}
where $\tau$ and $c_{r}$ are positive constant. Finally, (\ref
{zeta+1}) is deduced from
\[
\tilde{\mu}_{a,b} ( \exists x \in\tilde{V}_{n} , \tilde
{I}_{x}^{+} = 0 ) \leq\exp\biggl( - R(n)^{d} W_{n}(\eta)
\biggl( \tau c_{r} - \frac{d \ln(n/R(n))}{R(n)^{d} W_{n}(\eta)}
\biggr) \biggr) .
\]
\upqed\end{pf*}

This section ends with the proof of Property \ref{proper:markovian_tilde}.
\begin{pf*}{Proof of Property \protect\ref{proper:markovian_tilde}}
First, note that any event $A\in\mathcaligr{F}(\tilde{U})$ can be
written as a disjoint union of configurations of $\{-,+\}^{\tilde
{U}}$. So, it suffices to prove the identities (\ref{mu_tilde_mu}) and
(\ref{markovian_tilde}) for $A=\{\zeta\}$, $\zeta\in\{-,+\}^{\tilde
{U}}$.

Let us pick such a configuration $\zeta\in\{-,+\}^{\tilde{U}}$. The
set $\{f^{-1}(\zeta'), \zeta'\in\{-,+\}^{\tilde{V}}\}$ is a $\pi
$-system which generates the $\sigma$-algebra $f^{-1}(\mathcaligr
{F}(\tilde{V}))$. Hence, it is enough to prove
%
%
\begin{equation}
\label{pi_system1}
\mathbb{E}_{a,b} \bigl[\tilde{\mu}_{a,b} (
\zeta
| \mathcaligr{F}(\tilde{V}) ) \circ f \mathbh{1}
_{f^{-1}(\zeta')} \bigr]= \mathbb{E}_{a,b} \bigl[
\mathbh{1}
_{f^{-1}(\zeta)} \mathbh{1}_{f^{-1}(\zeta')} \bigr]
\end{equation}
for all $\zeta'\in\{-,+\}^{\tilde{V}}$ and
%
%
\begin{equation}
\label{pi_system2}
\mathbb{E}_{a,b} [\tilde{\mu}_{a,b} (
\zeta
| \mathcaligr{F}(\tilde{V}) ) \circ f ]=
\mathbb{E}
_{a,b} \bigl[\mathbh{1}_{f^{-1}(\zeta)} \bigr]
\end{equation}
(see, e.g., \cite{Williams}, page 84). Let us start with relation (\ref
{pi_system1}). For a configuration $\zeta'$ belonging to $\{-,+\}
^{\tilde{V}}$,
\begin{eqnarray*}
\mathbb{E}_{a,b} \bigl[\tilde{\mu}_{a,b} (
\zeta
| \mathcaligr{F}(\tilde{V}) ) \circ f \mathbh{1}
_{f^{-1}(\zeta')} \bigr]& = & \mathbb{E}_{a,b} [
\tilde
{\mu}_{a,b} ( \zeta | \mathcaligr{F}(\tilde{V})
) \circ f \mathbh{1}_{\zeta'} \circ f ]\\
& = & \tilde{\mathbb{E}}_{a,b} [\tilde{\mu}_{a,b}
(
\zeta | \mathcaligr{F}(\tilde{V}) ) \mathbh{1}
_{\zeta'} ]\\
& = & \tilde{\mathbb{E}}_{a,b} [\mathbh{1}_{\zeta}
\mathbh{1}_{\zeta'}
]\\
& = & \mathbb{E}_{a,b} \bigl[\mathbh{1}_{f^{-1}(\zeta)}
\mathbh{1}
_{f^{-1}(\zeta')} \bigr] .
\end{eqnarray*}
Relation (\ref{pi_system2}) is treated in the same way:
\begin{eqnarray*}
\mathbb{E}_{a,b} [\tilde{\mu}_{a,b} (
\zeta
| \mathcaligr{F}(\tilde{V}) ) \circ f ]& = &
\tilde{\mathbb{E}}_{a,b} [\tilde{\mu}_{a,b} (
\zeta | \mathcaligr{F}(\tilde{V}) ) ]\\
& = & \tilde{\mathbb{E}}_{a,b} [\mathbh{1}_{\zeta}
]\\
& = & \mathbb{E}_{a,b} \bigl[\mathbh{1}_{f^{-1}(\zeta)}
\bigr] .
\end{eqnarray*}
Now, let us prove (\ref{markovian_tilde}) with $A=\{\zeta\}$, $\zeta
\in\{-,+\}^{\tilde{U}}$. The set $\{\zeta', \zeta'\in\{-,+\}
^{\delta\tilde{U}}\}$ is a $\pi$-system which generates the $\sigma
$-algebra $\mathcaligr{F}(\delta\tilde{U})$. So, since the random
variables $\tilde{\mu}_{a,b}(\zeta|\mathcaligr{F}(\tilde{V}))$ and
$\tilde{\mu}_{a,b}(\zeta|\mathcaligr{F}(\delta\tilde{U}))$ have the
same expectation (equal to $\tilde{\mu}_{a,b}(\zeta)$), it suffices
to prove that, for any $\zeta'\in\{-,+\}^{\delta\tilde{U}}$,
%
%
\begin{equation}
\label{pi_system3}
\tilde{\mathbb{E}}_{a,b} [\tilde{\mu}_{a,b} (
\zeta | \mathcaligr{F}(\tilde{V}) ) \mathbh
{1}_{\zeta'}
]= \tilde{\mathbb{E}}_{a,b} [\tilde{\mu}_{a,b}
( \zeta | \mathcaligr{F}(\delta\tilde{U})
) \mathbh{1}_{\zeta'} ] .
\end{equation}
Let $\zeta'$ be a configuration of $\{-,+\}^{\delta\tilde{U}}$.
First, relation (\ref{mu_tilde_mu}) allows us to express the above
expectations according to the measure $\mu_{a,b}$:
\begin{eqnarray*}
\tilde{\mathbb{E}}_{a,b} [\tilde{\mu}_{a,b} (
\zeta | \mathcaligr{F}(\tilde{V}) ) \mathbh
{1}_{\zeta'}
]& = & \mathbb{E}_{a,b} [\tilde{\mu
}_{a,b}
( \zeta | \mathcaligr{F}(\tilde{V}) ) \circ f
\mathbh{1}_{\zeta'} \circ f ]\\
& = & \mathbb{E}_{a,b} \bigl[\mu_{a,b} (
f^{-1}(\zeta
) | f^{-1}(\mathcaligr{F}(\tilde{V})) ) \mathbh{1}
_{f^{-1}(\zeta')} \bigr] .
\end{eqnarray*}
Thus, let us denote by $\Sigma_{1}$ and $\Sigma_{2}$ the following sets:
\[
\Sigma_{1} = \bigcup_{y \in\delta\tilde{U}} B\bigl(y,R(n)+r\bigr)
\quad\mbox{and}\quad \Sigma_{2} = \biggl( \bigcup_{y \in\tilde{V}} B\bigl(y,R(n)+r\bigr)
\biggr) \Bigm\backslash\Sigma_{1} .
\]
Since $f^{-1}(\zeta')$ is a subset of $\Sigma_{1}$, we can write
\begin{eqnarray*}
&&\tilde{\mathbb{E}}_{a,b} [\tilde{\mu}_{a,b}
(\zeta | \mathcaligr{F}(\tilde{V}) ) \mathbh{1}
_{\zeta'} ]\\
&&\qquad = \mathop{\mathop{\sum}_{\omega\in\{-,+\}^{\Sigma_{1}}}}_
{f(\omega)=\zeta'} \mathbb{E}_{a,b} [\mu_{a,b}
( f^{-1}(\zeta) | f^{-1}(\mathcaligr{F}(\tilde
{V})) ) \mathbh{1}_{\omega} ]\\
&&\qquad = \mathop{\mathop{\sum}_{\omega\in\{-,+\}^{\Sigma_{1}}}}_
{f(\omega)=\zeta'} \sum_{\omega' \in\{-,+\}^{\Sigma
_{2}}} \mathbb{E}_{a,b} [\mu_{a,b} (
f^{-1}(\zeta
) | f^{-1}(\mathcaligr{F}(\tilde{V})) ) \mathbh
{1}_{\omega
\omega'} ] ,
\end{eqnarray*}
where $\omega\omega'$ is a configuration of $\{-,+\}^{\Sigma_{1}\cup
\Sigma_{2}}$. Now, the random variable\break $\mu_{a,b}(f^{-1}(\zeta)
|f^{-1}(\mathcaligr{F}(\tilde{V})))$ only depends on the vertices of
$\Sigma_{1}\cup\Sigma_{2}$. Hence,
\[
\mathbb{E}_{a,b} [\mu_{a,b} ( f^{-1}(\zeta
)
| f^{-1}(\mathcaligr{F}(\tilde{V})) ) \mathbh
{1}_{\omega
\omega'} ]= \mu_{a,b} ( f^{-1}(\zeta)
| \omega\omega' ) \mu_{a,b} ( \omega\omega'
) .
\]
Furthermore, the configurations belonging to $f^{-1}(\zeta)$ only
depend on the vertices of balls $B(y,R(n)+r)$, $y\in\tilde{U}$, and
by construction of the graph $\tilde{G}_{n}$ (and the use of the
$L_{\infty}$ norm) the following inclusion holds:
\[
\delta\biggl( \bigcup_{y \in\tilde{U}} B\bigl(y,R(n)+r\bigr) \biggr) \subset
\Sigma_{1} .
\]
So, the Markovian character of the measure $\mu_{a,b}$ applies. The
conditional probability $\mu_{a,b}(f^{-1}(\zeta) |\omega\omega')$
can be reduced to $\mu_{a,b}(f^{-1}(\zeta) |\omega)$ (see \cite
{Malyshev}, page 7). Combining the previous equalities, it follows that
\begin{eqnarray*}
&&\tilde{\mathbb{E}}_{a,b} [\tilde{\mu}_{a,b}
(\zeta | \mathcaligr{F}(\tilde{V}) ) \mathbh{1}
_{\zeta'} ]\\
&&\qquad = \mathop{\mathop{\sum}_{\omega\in\{-,+\}^{\Sigma_{1}}}}_
{f(\omega)=\zeta'} \mu_{a,b}
(f^{-1}(\zeta) | \omega) \sum_{\omega' \in\{-,+\}
^{\Sigma_{2}}} \mu_{a,b} ( \omega\omega' ) \\
&&\qquad = \mathop{\mathop{\sum}_{\omega\in\{-,+\}^{\Sigma_{1}}}}_
{f(\omega)=\zeta'} \mu_{a,b} (
f^{-1}(\zeta) | \omega) \mu_{a,b} ( \omega
) \\
&&\qquad = \mathop{\mathop{\sum}_{\omega\in\{-,+\}^{\Sigma_{1}}}}_
{f(\omega)=\zeta'} \mathbb{E}_{a,b} [\mu_{a,b}
( f^{-1}(\zeta) | f^{-1}(\mathcaligr{F}(\delta
\tilde{U})) ) \mathbh{1}_{\omega} ]\\
&&\qquad = \mathbb{E}_{a,b} \bigl[\mu_{a,b} (
f^{-1}(\zeta
) | f^{-1}(\mathcaligr{F}(\delta\tilde{U})) )
\mathbh{1}
_{f^{-1}(\zeta')} \bigr] .
\end{eqnarray*}
Finally, using a second time the relation (\ref{mu_tilde_mu}), we get
the desired identity:
\begin{eqnarray*}
\tilde{\mathbb{E}}_{a,b} [\tilde{\mu}_{a,b} (
\zeta | \mathcaligr{F}(\tilde{V}) ) \mathbh
{1}_{\zeta'}
]& = & \mathbb{E}_{a,b} [\tilde{\mu
}_{a,b}
( \zeta | \mathcaligr{F}(\delta\tilde{U}) )
\circ f \mathbh{1}_{\zeta'} \circ f ]\\
& = & \tilde{\mathbb{E}}_{a,b} [\tilde{\mu}_{a,b}
(\zeta | \mathcaligr{F}(\delta\tilde{U}) )
\mathbh{1}_{\zeta'} ] .
\end{eqnarray*}
\upqed\end{pf*}

%

%
\printaddresses

\end{document}